\input amstex
\magnification=1200
\documentstyle{amsppt}
\nologo
\baselineskip 20pt

\topmatter
\title On the Complete Integrability of Completely
Integrable Systems\endtitle 
\rightheadtext{}
\author Richard Beals and D. H. Sattinger \endauthor
\affil Yale University and University of Minnesota \endaffil
\leftheadtext{}
\footnote[]{Research
supported by NSF grants DMS-8916968 and DMS 8901607}

\abstract {The question of complete
integrability of evolution equations associated to  $n\times n$
first order isospectral operators is
investigated using the inverse scattering method.  It is shown that
for  $n>2$, e.g.  for the three-wave interaction, additional (nonlinear)
pointwise flows are necessary for the assertion of complete
integrability.  Their existence is demonstrated by constructing
action-angle variables.  This construction depends on the
analysis of a natural $2$-form and symplectic foliation for the
groups  $GL(n)$  and $SU(n)$.}
 \endabstract

\endtopmatter
\document
\baselineskip 15pt

\centerline{Comm. Math. Phys. {\bf 138}, 1991, 409-436}

\head 1.  Introduction \endhead

A classical Hamiltonian flow with  $2N$  degrees of freedom is said
to be completely integrable if it has $N$ independent integrals of the
motion which are in involution.  More generally,  $k$  independent
commuting Hamiltonian flows in a $2N$-dimensional manifold are said to
be a completely integrable family if there are  $N-k$  independent
integrals of the motions which are in involution, or equivalently
the  $N-k$  flows may be enlarged to a set of  $N$  independent commuting
flows.  By a theorem of Jacobi and Liouville, there then exist (at
least locally in phase space) a new set of canonical variables,
called action-angle variables, in which the flows are particularly
simple; set \cite{A} for a precise global version due to Arnold.
     
In recent years a number of nonlinear evolution equations,
beginning with the KdV equation, have been shown to have Hamiltonian
form on appropriate infinite-dimensional manifolds and to have an
infinite family of integrals of the motion which are in involution.
Such equations are commonly referred to as ``completely integrable'',
although it no longer makes sense to count half the number of
dimensions.  Nevertheless the inverse scattering 
method makes it possible to give a precise form to the
question of complete integrability and, indeed, to reduce it to a
question in a finite dimensional space.

Such results are known for the KdV hierarchy and the nonlinear
Schr\"odinger hierarchy; \cite{Ga}, \cite{ZF}.  The
scattering map both linearizes and decouples these flows, and action
angle variables have been obtained \cite{ZF}, \cite{ZM}. In both
these cases the pointwise 
dimension of the scattering data is  $2$, and no new pointwise flows
are needed to get to half the dimension.  Manakov \cite{Ma} obtained
action-angle variables for the $3$-wave interaction equation.

The inverse scattering method is based on analysis of an
associated linear spectral problem, and the associated hierarchy of
flows are isospectral for the linear operator.  For KdV the linear
operator is the  $1$-dimensional Schr\"odinger operator and for NLS
it is the  $2\times 2$  AKNS-ZS operator
$$
\frac{d}{dx}-zJ - q(x);\tag1.1
$$
here  $z$  is the spectral parameter, and  $J$  is a given constant
matrix.  In this paper we consider the  $n\times n$  version of
(1.1), under the assumption that $J$  is semisimple with distinct
eigenvalues and that  $q$  takes values in the range of  $ad~J$  and
vanishes rapidly at $\infty$.  The  three-wave interaction is an
example of an associated isospectral flow, when  $n=3$; \cite{ZM},
\cite{Ma}, \cite{Ka}.   The scattering and inverse scattering theory
for (1.1) has been considered by a number of authors: see \cite{ZS1},
\cite{AKNS} for  $n=2$  and \cite{Ma}, \cite{Ka}, \cite{Sh},
\cite{Ne}, \cite{Ge}, \cite{BY}, \cite{BC1}, \cite{Ca}. 

Each traceless matrix  $\mu$  which commutes with  $J$  generates a
hierarchy of isospectral flows of (1.1).  On the scattering side these
form an  $(n-1)$-parameter family of commuting pointwise fows.  The
pointwise dimension of the space of scattering data for (1.1) is
$n^2-n$.  We show that there is an appropriate Hamiltonian structure on
this space of pointwise data, and that the  $(n-1)$-parameter family is
completely integrable in the classical sense: it is part of an
$(n^2-n)/2$-parameter family of commuting Hamiltonian flows.  Both the
existence of an appropriate pointwise Hamiltonian structure and
complete integrability follow from a construction of Darboux
coordinates (coordinates which diagonalize the $2$-form) which are
action-angle variables for the flows of the hierarchies.  It should be
noted that the additional commuting point-wise flows needed when  $n$
is greater than $2$  are not linear on scattering data.
     
The results just described are obtained in the category of
complex manifolds and Hamiltonian structures.  We are also interested
in the real case.  The three-wave interaction, for example, is
associated to the operator (1.1) with  $J+J^{\ast}=0$  and the constraint
$q+q^{\ast}=0$;  on the scattering side the appropriate group is  $SU(3)$
rather than  $SL(3,C)$.  We show that one can find real Darboux
coordinates for scattering data to provide action-angle variables for
the $3$-wave interaction and the other flows of the hierarchies.  The
canonical transformation to action-angle variables is not algebraic in
this case: it requires the Liouville method and elliptic functions.
Manakov \cite{Ma} used a different method to obtain action-angle
variables for the $3$-wave interaction which have a simple form but
which are nonlocal functions of the entries of the scattering matrix
$s$  of section 3; the associated flows are also nonlocal in $s$.

Our analysis of the Hamiltonian structures (symplectic form,
Poisson brackets) leads to a natural closed $2$-form of rank  $n^2-n$
on  $GL(n)$, and a natural symplectic foliation of  $GL(n)$.  
The reduction  $J+J^{\ast}=0$, $q+q^{\ast}=0$  leads to consideration
of  $SU(n)$ in place of  $GL(n)$. The induced Poisson bracket is not
a Poisson-Lie structure \cite{Dr}, since it is not degenerate at the
identity element.  However it was pointed out to us by Lu \cite{Lu}
that our structure is the translate by a Weyl element of a
Poisson-Lie structure which is the classical limit of a quantum group
structure described by Drinfeld \cite{Dr}.

The plan of the paper is the following.  In section 2 we review
the Hamiltonian structure and hierarchy of flows associated to the
operator (1.1).  The scattering theory for the case  $J+J^{\ast}=0$  is
reviewed in section 3.  We then compute the Poisson bracket for
scattering data and state the main results on existence of Darboux
coordinates and complete integrability.  In section 4 we introduce and
analyze the $2$-form on  $GL(n)$  and obtain Darboux coordinates.  A
symplectic foliation of  $GL(n)$ is introduced in section 4, and we
calculate the associated Poisson bracket and the Hamiltonians for a
family of linear flows.
    
The algebraic results of sections 4 and 5 are used in section 6 to
prove the results on Darboux coordinates and complete integrability
for scattering data which were stated in section 3. The case of $SU(3)$
is taken up in section 7; complete integrability of the three-wave
interaction is a consequence.  In section 8 we show that the results
stated in section 2 remain valid valid without the restriction  
$J+J^{\ast}=0$.

We are grateful for discussions with V. Zurkowski, J. H. Lu, G.
Moore, and G. Zuckerman.  B. Konopelchenko drew our attention to the
work of Manakov.

\head 2.  Symplectic structure of Hamiltonian hierarchies \endhead

We consider Hamiltonian hierarchies of flows associated to the
first order differential operator
$$
\frac{d}{dx} - zJ - q(x),\quad z\in \bold C,\tag2.1
$$
where  $J$  is a constant  $n\times n$  semisimple matrix;  $q(x)$
is an  $n\times n$  matrix whose entries  $q_{jk}$  belong to the
Schwartz class  $\Cal S(\bold R)$; and, for each  $x$,  $q(x)$  lies
in the range of  $ad J$.  We denote by  $P$
the linear space of all such  $q$; thus  $P = \Cal S(\bold R; ad
J(M_n))$, where  $M_n = M_n(\bold C)$  is the space of  $n\times n$
matrices, with the Schwartz topology.  We use the following inner
product on  $P$:
$$
\langle q,p\rangle = \int_{\bold R} tr[q(x)p(x)]dx.\tag2.2
$$
Since  $P$  is a linear space we may identify it with its tangent space.
We denote tangent vectors (at a given point  $q$) by  $\dot q$.  Associated
with (2.1) is the closed $2$-form
$$
\Omega_P = \frac12 \int_{\bold R} tr[\delta q(x) \wedge [ad
J]^{-1}\delta q(x)]dx \tag2.3
$$
where $\delta q(x)$  denotes the linear functional taking  $\dot q$
to  $\dot q(x)$  and  $[adJ]^{-1}$  maps to the range of  $adJ$,  on
which  $adJ$  is injective.  Thus
$$
\Omega_P(\dot q_1, \dot q_2 ) = \frac12\int_{\bold R} tr[\dot
q_1(x)[ad J]^{-1}\dot q_2 (x) - \dot q_2(x)[ad J]^{-1}\dot
q_1(x)]dx \tag2.3$^\prime$
$$
Since the inner product is non-singular,  $\Omega_P$  is symplectic.  Note
that when  $J+J^{\ast} = 0$, and we restrict to the set $\{q\in P:q+q^{\ast} =
0\}$, then the form  $\Omega_P$  is real.

We shall work with the case in which  $J$  is diagonal with
distinct eigenvalues:  $J = \operatorname{diag}(i\lambda_1, i\lambda_2,\dots
i\lambda_n)$ .  In this case
$$
\Omega_{P} = \int_{\bold
R}\sum_{j<k}\frac{1}{i(\lambda_k-\lambda_j)}\delta
q_{jk}(x)\wedge\delta q_{kj}(x)dx.\tag2.4
$$

A Poisson bracket is associated to the symplectic form  $\Omega_P$  in
the standard way.  If  $F$  is a functional on  $P$  which is Frechet
differentiable, and  $\dot q$  is a tangent vector, we write
$$
[\dot q F](q) = \frac{d}{d\varepsilon}\bigm|_{\varepsilon =0} F(q +
\varepsilon\dot q) = \langle\frac{\delta F}{\delta q}(q),\dot q
\rangle,\tag2.5 
$$
i.e. we identify  $\delta F/\delta q$  with the gradient of  $F$.
The Hamiltonian vector field associated to  $F$ , denoted   $H_F$, is
then defined by
$$
\Omega_P(H_F,\dot q) = - \dot q F .\tag2.6
$$
This definition and (2.3$'$) imply that  $H_F  = [J,\delta F/\delta
q]$.  The Poisson bracket \cite{Ne} is then given by 
$$
\{F,G\}_P = H_FG = -\Omega_P(H_G,H_F) = \int_{\bold R}tr\left([J,\frac{\delta
F}{\delta q}]\frac{\delta G}{\delta q}\right)dx. \tag2.7
$$

There is an  $(n-1)$-parameter family of hierarchies of commuting
Hamiltonian flows in  $P$, defined as follows.  Let  $\mu$  be a constant
matrix with
$$
tr\mu = 0,\quad [J,\mu] = 0 ,
$$
and associate to  $q$  in  $P$  a sequence of matrix-valued functions
$F_k$  defined recursively by
$$\align
F_0(x) &= \mu;\\
[J,F_{k+1}] &=  \frac{dF_k}{dx} + [q, F_k];\quad\lim\limits_{x\to
-\infty} F_{k+1}(x) = 0.
\endalign
$$
The  $F_k$  depend nonlinearly on  $q$  for  $k>1$
$(k>2,\quad\text{if}\quad n = 2)$.  Various formal and rigorous
versions of the following are well-known.

\proclaim{Theorem 2.1}  \cite{\bf Sa, \bf BC2\rm, \bf BC3\rm}.  Each
$F_k(q)$  is a  polynomial in  $q$  and its derivatives of order less
than  $k$.  The hierarchy of flows defined by
$$
\dot q = [J, F_{k+1}(q)]\tag2.8
$$
are Hamiltonian with respect to  $\Omega_P$  and the Hamiltonians
are in involution with respect to the Poisson bracket  $\{~,~\}_P$.
\endproclaim

We shall discuss the Hamiltonians for these flows later.

It is also well-known that the scattering transform linearizes
the flows (2.8).  We discuss this in the next section.

\head 3.  The scattering transform; symplectic structure on scattering
data \endhead

We summarize here the basic results of scattering theory for the
operator (2.1); cf \cite{BC1}.  In this section we assume
$$
J = \operatorname{diag}(i\lambda_1,\dots,i\lambda_n),\quad\lambda_j\in\bold R,
\quad\lambda_1 > \lambda_2 > \dots > \lambda_n. \tag3.1
$$
For a given  $q$  in  $P$,   we seek a matrix-valued solution of the
spectral problem
$$
\frac{\partial}{\partial x}\psi(x,z) = zJ\psi(x,z) + q(x) \psi(x,z),
\quad z\in\bold C.\tag3.2
$$
which is normalized by the asymptotic conditions
$$
\lim\limits_{x\to -\infty}\psi(x,z)\exp(-xzJ) = 1,\quad\underset
x\rightarrow +\infty\to{\lim\sup}\Vert\psi(x,z)exp(-xzJ)\Vert
<\infty. \tag3.3
$$
If  $\int\Vert q(x)\Vert dx < 1$,  then there is a unique solution to
(3.2), (3.3), and it has a limit
$$
\lim\limits_{x\to +\infty}\exp(-x\xi J) \psi(x,\xi) = s(\xi),
\quad\xi\in\bold R .\tag3.4
$$
The transformation  $q\mapsto s=s(\cdot;q)$  is one of two versions of the
{\it scattering transform}, and  $s$  is called the
{\it scattering matrix}.

Still assuming   $\int\Vert q(x)\Vert dx<1$,  the solutions
$\psi(x,z)$  for non-real  $z$  are holomorphic and have limits on
$\bold R$  which are related by
$$
\psi(x,\xi + i0) = \psi(x,\xi-i0)v(\xi),\quad\xi\in\bold R.\tag3.5
$$
To describe the target spaces for the maps  $q\mapsto s$  and
$q\mapsto v$,  we define the spaces
$$\align
SL^{\pm} &= \{a\in SL(n,C) : a_{jk} = 0
\quad\text{if}\quad\pm(j-k)>0\};\tag3.6\\
SL^{\pm}_0 &= \{a\in SL^{\pm} : \operatorname{diag}(a) = 1\};\\
SL_{\ast} &= (SL^+\cdot SL^-)\cap (SL^-\cdot SL^+).
\endalign
$$
This means that  $SL_{\ast}$  consists precisely of those  $s$  in  $SL(n) =
SL(n,\bold C)$  which have two (unique) triangular factorizations
$$
s = s_+ v^{-1}_+ = s_- v^{-1}_-,\quad s_{\pm}\in SL^{\pm},\quad v_{\pm}\in
SL^{\mp}_0.\tag3.7
$$

\proclaim{(3.8)  Definition}  $SD$  is the set of matrix-valued
functions $s:\bold R\to SL(n)$  with the properties
\roster
\item"(3.8a)"   $s$  is smooth and bounded; each derivative has an
asymptotic expansion in powers of  $\xi^{-1}$  as  $\vert\xi\vert\to\infty$;
\item"(3.8b)"   $s$  takes values in  $SL_{\ast}$, so it factors as
$s(\xi) = s_{\pm}(\xi)v_{\pm}(\xi)^{-1}$;
\item"(3.8c)"   the diagonal-matrix-valued functions
$\delta_{\pm}(\xi) = \operatorname{diag} s_{\pm}(\xi)$  are
the boundary values of a diagonal-matrix-valued function
which is bounded, holomorphic, and invertible in  $\bold
C\backslash\bold R$.
\endroster
\endproclaim

\proclaim{(3.9)  Definition}  $SD'$  is the set of pairs of matrix-valued
functions  $(v_+,v_-), v_{\pm} :\bold R\to SL^{\mp}$, with the properties
\roster
\item"(3.9a)"   each entry of  $v_{\pm}-1$  belongs to  $\Cal S(\bold R)$;
\item"(3.9b)"   the upper principal minors of the matrix-valued function
$v=v^{-1}_-v_+$  are non-zero and have winding number zero.
\endroster
\endproclaim

Condition (3.9b) excludes discrete scattering data (bound
states); throughout this paper we consider only potentials with
purely continuous scattering data.

We equip  $SD$  and  $SD'$  with the Schwartz topologies.

\proclaim{Theorem 3.1} (\cite{Sh}, \cite{B-Y}, \cite{BC1}).  The map
$q\mapsto s$ 
is a diffeomorphism from a neighborhood of  $0$  in  $P$  onto a
neighborhood of  $1$  in  $SD$, and it extends to map an open set in
$P$  bijectively to a dense open set in  $SD$.

The matrix function  $v$  in (3.5) is related to the matrix
function  $s$  by the factorizations (3.7); in fact
$$
v=v^{-1}_-v_+=s^{-1}_-s_+,\tag3.10
$$
and this equation uniquely determines  $v_{\pm}$   from  $v$.  The map  $q
\mapsto (v_+,v_-)$  is a diffeomorphism from a neighborhood of the
origin in  $P$ 
onto a neighborhood of  $(1,1)$  in  $SD'$, and it extends to map an open
set in  $P$  bijectively onto a dense open set in  $SD'$.
\endproclaim

\demo{(3.11)  Remark}  It is nearly implicit that  $SD$  and  $SD'$  are
diffeomorphic.  The factorizations (3.7) and (3.10) determine  $v$  from
$s$.  Conversely, write  $s_{\pm}=\delta_{\pm}t_{\pm}$  with
$\delta_{\pm}$  diagonal and  $t_{\pm}$  in  $SL^{\pm}_0$.
The factorization (3.10) gives  $v^{-1}_-v_+ =
t^{-1}_-(\delta^{-1}_-\delta_+)t_+$, showing that  $t_{\pm}$  and
$\delta^{-1}_-\delta_+$  are determined algebraically from  $v$.  The
holomorphy  properties (3.9)  show that the factors
$\delta_-,\delta_+$  can be obtained from  $\delta^{-1}_-\delta_+$
by solving a Riemann-Hilbert factorization problem, and thus
$s_{\pm}$   and  $s$  itself can be obtained from  $v$  or from the
pair  $(v_+,v_-)$.
\enddemo

\proclaim{Proposition 3.2} (\cite{G}, \cite{BC1}).  The pull-back of the
$2$-form  $\Omega_P$  of (2.3) under the inverse of the scattering
transform is
$$
\Omega_S = \frac{1}{4\pi i}\int_{\bold R} tr[v_-(\delta
v)v^{-1}_+\wedge s^{-1}\delta s].\tag3.12
$$
\endproclaim

It will be convenient to have a somewhat different formulation.

\proclaim{Proposition 3.3}  The $2$-form  $\Omega_S$  can be written
$$
\Omega_S = \frac{1}{4\pi i}\int_{\bold R} tr[v^{-1}_+\delta v_+\wedge
s^{-1}_+\delta s_+ - v^{-1}_-\delta v_-\wedge
s^{-1}_-\delta s_-].\tag3.13
$$
\endproclaim

\demo{Proof}  Since  $v=v^{-1}_-v_+$  and  $s=s_{\pm}v^{-1}_{\pm}$,
it follows that
$$\align
v_-(\delta v)v^{-1}_+ &= (\delta v_+)v^{-1}_+ - (\delta
v_-)v^{-1}_-;\\
s^{-1}\delta s &= s^{-1}[\delta s_{\pm}v^{-1}_{\pm} -
s_{\pm}v^{-1}_{\pm}\delta v_{\pm}v^{-1}_{\pm}] =
v_{\pm}s^{-1}_{\pm}(\delta s_{\pm})v^{-1}_{\pm} - (\delta
v_{\pm})v^{-1}_{\pm}.
\endalign
$$
Then
$$
tr[(\delta v_{\pm})v^{-1}_{\pm}\wedge s^{-1}\delta s] =
tr[v^{-1}_{\pm}\delta v_{\pm}\wedge s^{-1}_{\pm}\delta s_{\pm}],
$$
since  $(\delta v_{\pm})v^{-1}_{\pm}$  is strictly upper or lower
triangular.
\enddemo

Next we consider the image under the scattering transformation of
the Poisson bracket  $\{~,~\}_P$   of (2.7); equivalently, this is the
Poisson bracket associated with the $2$-form  $\Omega_S$  on scattering data.
As usual, we may consider the entries of the
scattering matrix  $s(\xi) = s(\xi;q)$  to be functionals on  $P$  and
compute the corresponding bracket
$$
\{s_{jk}(\xi),s_{\ell m}(\eta)\}_{_{S}}(s) = \{s_{jk}(\xi),s_{\ell
m}(\eta)\}_{_{P}}(q),\quad\xi,\eta\in\bold R.
$$
There are two problems here.  First, the gradients  $\delta
s_{jk}(\xi)/\delta q$  do not decay, so the formula (2.7) does not
have an absolutely convergent integrand and it is necessary to use a
regularization such as
$$
\lim\limits_{N\to\infty}\int^N_{-N}tr\left([J,\frac{\delta
s_{jk}}{\delta q}(\xi)]\frac{\delta s_{\ell m}}{\delta q}(\eta)\right)dx.
$$
Second, even this limit exists only in the sense of distributions in
the two variables  $\xi,\eta$.  Thus the precise meaning of the calculation
is this: for any pair of test functions  $u, w$  in
$C^{\infty}(\bold R)$  one considers the pair of functionals
$$
F(q) = \int_{\bold R} s_{jk}(\xi)u(\xi)d\xi,\quad G(q) = \int_{\bold
R}s_{\ell m}(\xi)w(\xi)d\xi,\quad s=s(\cdot;q).
$$
Then formally one has
$$\align
\{F,G\} &= \lim\limits_{N\to\infty} \int^N_{-N}tr(
[J,\frac{\delta F}{\delta q}]\frac{\delta G}{\delta q})dx \tag3.14\\
&= \lim\limits_{N\to\infty}\iiint^N_{-N} tr\left([J,\frac{\delta s_{jk}(\xi)}
{\delta q}]\frac{\delta s_{\ell m}(\xi)}{\delta q}\right)
u(\xi)w(\eta)dxd\xi d\eta\\  
&= \iint\{s_{jk}(\xi),s_{\ell m}(\eta)\}_S u(\xi)w(\eta)d\xi d\eta
\endalign
$$
as the defining equation for the distribution
$\{s_{jk}(\xi),s_{\ell m}(\eta)\}\in\Cal D'(\bold R\times\bold R)$.
The following calculation is standard; see \cite{Ma} for the
$3\times 3$ case and  \cite{Sk}, \cite{KD} for  $R$-matrix
formulations.

\proclaim{Proposition 3.4}  The distribution defined by (3.14) is given
explicitly by
$$\align
\{s_{jk}(\xi),s_{\ell m}(\eta)\} &= \pi i s_{jm}(\xi)s_{\ell
k}(\eta)[sgn(\ell -j)-sgn(m-k)]\delta(\xi-\eta)\tag3.15\\
&+ s_{jk}(\xi)s_{\ell m}(\eta)[\delta_{j\ell} - \delta_{km}]p.v.
\frac{1}{\xi-\eta}
\endalign
$$
{\it where we take}  $sgn(0) = 0$  and  $p.v.$  {\it denotes the
principal value}.
\endproclaim

\demo{Proof}  The variation of  $s$  with respect to  $q$  is
$$
\dot s(\xi) = \int_{\bold R}s(\xi)\psi(x,\xi)^{-1}\dot
q(x)\psi(x,\xi)dx;\tag3.16
$$
\cite{BC3, (2.45)}.  Here the  $\psi$  are the eigenfunctions (3.2),
normalized at  $x=-\infty$.  We write  $\tilde{\psi}(x,\xi) =
\psi(x,\xi)s(\xi)^{-1}$,  which
is normalized at  $x=+\infty$.  With  $F$  as above, an easy
calculation using (2.2), (2.5), and (3.16) shows that
$$
\frac{\delta F}{\delta q}(x) = \int_{\bold
R}\psi(x,\xi)e_{kj}\tilde{\psi}(x,\xi)^{-1}u(\xi)d\xi.
$$
A similar formula holds for  $G$,  so (3.14) becomes
$$\multline
\{F,G\}_P = \lim\limits_{N\to\infty}\int^N_{-N}tr([J,\frac{\delta F}{\delta
q}]\frac{\delta G}{\delta q})dx =\\
\lim\limits_{N\to\infty}\iiint^N_{-N}tr([J,\psi(x,\xi)e_{kj}
\tilde{\psi}(x,\xi)^{-1}]\psi(x,\eta)e_{m\ell}\tilde{\psi}
(x,\eta)^{-1})u(\xi)w(\eta)dx d\xi d\eta
\endmultline\tag3.17
$$
We use the identity
$$
\frac{1}{\xi-\eta}\frac{d}{dx}[\tilde{\psi}(x,\eta)^{-1}\psi(x,\xi)]
= \tilde{\psi}(x,\eta)^{-1} J\psi(x,\xi)
$$
and the properties of the trace to conclude from (3.17) that
$$\align
\{s_{jk}(\xi),s_{\ell m}(\eta)\}_S &=
\lim\limits_{N\to\infty}\frac{1}{\xi-\eta}
tr[e_{kj}g(\xi,\eta,N)e_{m\ell}g(\eta,\xi,N) 
\tag3.18\\
&- e_{kj}g(\xi,\eta,-N)e_{m\ell}g(\xi,\eta,-N)]
\endalign
$$
in the sense of distributions, where    $g(\xi,\eta,x) =
\tilde{\psi}(x,\xi)^{-1}\psi(x,\eta)$.
Now
$$\align
g(\xi,\eta,x) &\approx \exp[x(\eta-\xi)J]s(\eta)\quad\text{as}\quad
x\to +\infty;\\
&\approx s(\xi)\exp[x(\eta-\xi)J]\quad\text{as}\quad x\to -\infty.
\endalign
$$
Thus the right side of (3.18) is
$$
\lim\limits_{N\to\infty}\frac{1}{\xi-\eta}[s_{jm}(\eta)s_{\ell
k}(\xi)e^{iN(\xi-\eta)(\lambda_{\ell}-\lambda_{j})}-
s_{jm}(\xi)s_{\ell k}(\eta)e^{iN(\xi-\eta)(\lambda_{m}-\lambda_{k})}]
\tag3.19
$$
There is no singularity in (3.19) since the term in brackets vanishes
at $\xi=\eta$;  therefore we may replace the expression in (3.19) by the
principal value integral, i.e. letting the distribution act as the
limit as  $\varepsilon\downarrow 0$ of the integral over the region
$\vert\xi-\eta\vert >\varepsilon$.  This allows us to decouple the
two terms in (3.19) and use the identity
$$
\lim\limits_{N\to\infty} p.v. \frac{1}{\xi-\eta} e^{iaN(\xi-\eta)} =
\pi i sgn(a)\delta(\xi-\eta)\quad\text{if}\quad a\in\bold R\backslash
0.
$$
to deduce (3.15) from (3.19).
\enddemo

As is well-known, the flows of Theorem 2.1 become linear on the
scattering side.

\proclaim{Proposition 3.5}  The potential  $q(\cdot,t)$  evolves
according to (2.8) if and only if the scattering data evolve
according to
$$
\frac{\partial}{\partial t}s(\xi,t) = \xi^k[\mu,s(\xi,t)],\quad
\frac{\partial}{\partial t}v_{\pm}(\xi,t) =
\xi^k[\mu,v_{\pm}(\xi,t)].\tag3.20
$$
\endproclaim

For a proof, see for example \cite{BC2}, \cite{BC3}.  These flows are
Hamiltonian (with the same Hamiltonian functions as in the original
variables) with respect to the 
symplectic form  $\Omega_S$  on scattering data since the structure has
simply been transferred from  $P$  to  $SD$  or $SD'$.  We wish to emphasize
that on the scattering side the flows (2.8) are not only linearized
but decoupled for different values of $\xi,\eta$;  equivalently, the
Hamiltonian vector fields act in a pointwise fashion on the entries of
$s$  or of $(v_+,v_-)$.  This allows us to reduce the question of complete
integrability of the family of flows (2.8)/(3.20) to a
finite-dimensional problem.

This integrability question is related to certain problems and
questions concerning the symplectic and Poisson structures on
scattering data.  Observe that the $2$-form  $\Omega_S$  lifts to the
loop space
$$
L(SL_{\ast}) = \{a:\bold R\to SL_{\ast} : \text{entries of }
a-1 \text{ belong to } \Cal S(\bold R)\}.
$$
More precisely  $\Omega_S$ is the pullback to  $SD$  of the $2$-form
defined by (3.12) or (3.13) on  $L(SL_{\ast})$.  Moreover, these are
pointwise formulas, in the obvious sense: they express the form as
the direct integral of forms computed pointwise from entries of  a
in $L(SL_{\ast})$.  However, the 
form is not symplectic on  $L(SL_{\ast})$.  (This can be seen from the fact
that the dimension of $SL$  is  $n^2-1$  but the rank of the pointwise form
in (3.13), as we show in the next section, is  $n^2-n$.)

The alternative space of scattering data, $SD'$, is itself a loop
space, whose pointwise dimension is  $n^2-n$, and the form  $\Omega_S$  is
symplectic on  $SD'$; however the computation of  $s$  from
$(v_+,v_-)$ involves the solution of a Riemann-Hilbert problem, so
that (3.12) does not express  $\Omega_S$  as a direct integral of
pointwise $2$-forms on the fiber of the loop space; in fact it has no
such pointwise expression when  $n$  is larger than $2$.

The Poisson bracket  $\{~,~\}_S$   lifts to the whole loop group
$L(SL(n))$  using the same formula (3.15), but the nonlocal term
involving  $p.v.(\xi-\eta)^{-1}$ shows that  $\{~,~\}_S$   is also
not a direct integral of a pointwise defined bracket, unlike the
bracket  $\{~,~\}_P$  on the space  $P$  of potentials.  This
corresponds to the fact that the submanifold  $SD\subset L(SL(n))$
is determined in part by the nonlocal constraint (3.8c).  The
nonlocal term does not vanish when  $\{~,~\}_S$
is considered as a bracket on  $SD'$ if  $n > 2$, which proves the
earlier assertion that  $\Omega_S$  has no pointwise expression on  $SD'$.

The following summarizes these observations.

\proclaim{Proposition 3.6}  $\Omega_S$  is a closed, pointwise $2$-form
but not a symplectic $2$-form on the loop space  $L(SL_{\ast})$.
$\Omega_S$
is a symplectic $2$-form but not a pointwise $2$-form on  $SD'$.

$\{~,~\}_S$  is not a pointwise Poisson bracket on  $L(SL(n))$, nor
on  $SD$, nor on  $SD'$ (when  $n > 2$).
\endproclaim

Our main positive results involve choosing new coordinates to
overcome the limitations in Proposition 3.6.  To state them we must
extend the notion of a distribution-valued Poisson bracket beyond the
coordinate functions themselves.  Suppose  $f, g$  are in
$C^{\infty}(SL(n))$ and  $u,w$  are test functions.  We define
functionals on  $SD$, or on  $P$,  by
$$\alignat2
f_{\xi}(s) &= f(s(\xi)),  &&\qquad  F(s) = \int f(s(\xi))u(\xi)d\xi.\\
g_{\xi}(s) &=g(s(\xi)),  &&\qquad  G(s) = \int g(s(\xi))w(\xi)d\xi.
\endalignat
$$
Thus the coordinate functions  $a_{jk}(s) = s_{jk},\quad a_{\ell m}(s) =
s_{\ell m}$  give rise to the functionals   $s_{jk}(\xi), s_{\ell
m}(\eta), F$,  and  $G$  considered above.  Then as before we obtain
the distribution  $\{f_{\xi},g_{\eta}\}\in\Cal D'(\bold R\times\bold R)$
by the formal calculation
$$
(F,G)_P = \iint\{f_{\xi},g_{\eta}\}_s u(\xi)w(\eta)d\xi d\eta.
$$
This leads also to the expression
$$
\{f_{\xi},g_{\eta}\}_S = \sum\frac{\partial f}{\partial
s_{jk}}\frac{\partial g}{\partial s_{\ell m}}\{s_{jk}(\xi),s_{\ell
m}(\eta)\}_S.
$$

\proclaim{Theorem 3.7}  There are functions  $p_{\nu},q_{\nu},1\leq
\nu\leq \frac12(n^2-n)$ which are defined and holomorphic on a dense open
algebraic subset of $SL(n)$  and which have the following properties.
\roster
\item"(3.21)"  The map  $s\mapsto (p_1\circ s, p_2\circ
s,\dots,q_1\circ s, q_2\circ s,\dots)$ is an injection
from a dense open set in  $SD$  into  $C^{\infty}(\bold R;\bold
C^{n^{2}-n})$.
\item"(3.22)"  $\Omega_S=\sum_{\nu}\delta(p_{\nu}\circ s)\wedge
\delta(q_{\nu}\circ s)$
\item"(3.23)"  $\{p_{\mu,\xi},p_{\nu,\eta}\}_S =
\{q_{\nu,\xi},q_{\nu,\eta}\}_S = 0$;\quad $\{p_{\mu,\xi},q_{\nu,\eta}\}_S
= \delta_{\mu\nu}\delta(\xi-\eta)$.
\item"(3.24)"  For any  $f$  in  $C^{\ast}(SL(n))$, the distributions
$\{p_{\nu,\xi},f_{\eta}\}$  have support at  $\xi = \eta$,  all  $\nu$.
\endroster
\endproclaim

In other words, the functions  $p_{\nu}\circ s, q_{\nu}\circ s$   are
global Darboux coordinates on the manifold  $SD$  of scattering data.
The additional fact (3.24) implies that Hamiltonians which are
functions of the  $p_{\nu}$  give pointwise vector fields on SD.
These same functions  $p_{\nu}$  provide a strong positive answer to
our question about complete integrability in the  $SL(n)$ case, as follows.

\proclaim{Theorem 3.8}  The functions  $p_{\nu}$   of Theorem 3.2 may
be chosen so that for each traceless diagonal matrix  $\mu$,  the
Hamiltonian for the flow (2.8), (3.20) is a linear combination of the
functionals 
$$
\int_{\bold R}\xi^{k}p_{\nu}(s(\xi))d\xi.\tag3.25
$$
\endproclaim

Thus the family of flows (3.20), which is determined pointwise by
the $(n-1)$-parameter family of traceless diagonal matrices, is imbedded
in the $\frac12(n^2-n)$-parameter family of flows generated by the
$p_{\nu}$'s.  In fact the functions  $p_{\nu},q_{\nu}$  of Theorem
3.2 provide action-angle variables for the flows (3.20).

Theorems 3.7 and 3.8 are proved in section 6.

In this section we have made two restrictive assumptions about
$J$ :  that the eigenvalues are distinct and that they
lie on a line through the origin in  $\bold C$.  The same results
hold without the second assumption, as we show in section 8.

Another situation arises with {\it reductions}, i.e.the imposition of
restrictions on the potentials  $q$. The most interesting example in
the present context is the restriction  $q+q^{\ast} = 0$, still assuming
(3.1).  For such  $q$, the scattering data satisfies the corresponding
constraints
$$
s(\xi) \quad\text{belongs to}\quad SU(n);\qquad v_+(\xi)^{\ast}v_-(\xi) = 1.
$$
The pullbacks of the $2$-forms are still symplectic.  When  $n=3$  the
simplest associated nonlinear evolution equation is the $3$-wave
interaction.  (Note that one needs  $\mu +\mu^{\ast} =0$  in (2.8),
(3.20) to preserve the constraints.)   Here our manifolds are real,
and we need the functions  $p_{\nu},q_{\nu}$  of Theorem 3.2  to be
real in order to preserve the structure.

\proclaim{Theorem 3.9}  The functions  $p_{\nu},q_{\nu}$   of
Theorems 3.7, 3.8 can be chosen to be real on  $SU(3)$  in the case
$n= 3$.  In particular, the three-wave interaction is a completely
integrable Hamiltonian evolution in the strong sense.
\endproclaim

This result is proved in section 7.

\head  4.  A $2$-form in $GL(n)$ \endhead

     Let  $GL(n)$  denote either  $GL(n,\bold R)$  or  $GL(n,\bold
C)$; in the latter case the functions and forms to be considered are
complex-valued.  The key steps in deriving Theorems 3.2, 3.3, and 3.4
involve an analysis of a $2$-form in  $GL(n)$ : the integrand in
(3.13).  As in section 3 we introduce the matrix spaces
$$\align
GL^{\pm} &= \{a\in GL(n) : a_{jk} = 0 \quad\text{if}\quad \pm(j-k) > 0 \};\\
GL^{\pm}_0 &= \{a\in GL^{\pm} : \text{diag}(a) = 1\};\\
GL_{\ast} &= (GL^+\cdot GL^-)\cap (GL^-\cdot GL^+).
\endalign
$$
Again, $GL_{\ast}$  consists of those elements of  $GL(n)$  with
factorizations
$$
a=a_+v^{-1}_+ = a_-v^{-1}_-,\quad a_{\pm}\in GL^{\pm},\quad v_{\pm}\in
GL^{\mp}_0.\tag4.1
$$
We consider  $a_{\pm},v_{\pm}$  here as functions of  $a$  in  $GL(n)$; then
$da_{\pm}$  and $dv_{\pm}$  are matrices of  $1$-forms on  $GL(n)$
and we may define a $2$-form by
$$
\Omega = tr[v^{-1}_+dv_+\wedge a^{-1}_+da_+ - v^{-1}_-dv_-\wedge
a^{-1}_-da_-].\tag4.2
$$
At the identity,  $\Omega = 2\sum\limits_{j<k}da_{jk}\wedge da_{kj}$,
so  $\Omega$  has rank  $\geq n^2-n$  on an open set.

\proclaim{Theorem 4.1}  On a dense open algebraic subset of
$GL_{\ast}$  the $2$-form  $\Omega$ has a representation
$$
\Omega = \sum^N_{\nu =1}dp_{\nu}\wedge dq_{\nu},\qquad N=\frac12
(n^2-n),\tag4.3
$$
where  $p_{\nu}$  and  $q_{\nu}$  are analytic (holomorphic) and the
$1$-forms  $dp_{\nu}, dq_{\nu}$  are independent.  In particular,
$\Omega$  is closed and generically has rank  $n^2-n$.
\endproclaim

The proof of this theorem is given after Theorem 4.4.  The strategy
is to obtain the general result by a reduction to the case  $n=2$

\proclaim{Lemma 4.2}  For  $n=2$,
$$
\Omega = d\log \left[\frac{a_{11}a_{22}}{\triangle}\right]\wedge
d\log\left[\frac{a_{21}}{a_{12}}\right], 
\quad\triangle = \det a.\tag4.4
$$
\endproclaim

\demo{Proof}  If  $a$  is in  $GL_{\ast}(2)$  then
$$\align
v_+ &= \bmatrix 1 & 0\\-a_{21}/a_{22} &1\endbmatrix, a_+ =
\bmatrix \triangle/a_{22} &a_{12}\\ 0 & a_{22}\endbmatrix,\\
v_- &= \bmatrix 1 & -a_{12}/a_{11}\\ 0 & 1\endbmatrix, v_- =
\bmatrix a_{11} & 0\\ a_{21} & \triangle/a_{11}\endbmatrix
\endalign
$$
A direct calculation gives
$$\align
\Omega &= d\left[\frac{a_{21}}{a_{22}}\right]\wedge
\left[\frac{a_{12} da_{22}}{\triangle}
-\frac{a_{22}da_{12}}{\triangle}\right]\\
&-d\left[\frac{a_{12}}{a_{11}}\right]\wedge
\left[\frac{a_{21}da_{11}}{\triangle}
-\frac{a_{11}da_{21}}{\triangle}\right] \endalign
$$
and some further manipulation leads to (4.4).
\enddemo

\demo{(4.5)  Remark}  Given  $1\leq j<k\leq n$, consider the subset
of  $GL_{\ast}$  consisting of those  $a$  whose only non-vanishing
off-diagonal entries occur in the  $(j,k)$  and  $(k,j)$  places.
The pullback of $\Omega$ to this subset has, by an analogous
computation, the form
$$
\Omega = d\log\left[\frac{a_{jj}a_{kk}}{\Delta}\right]\wedge
d\log\left[\frac{a_{kj}}{a_{jk}}\right], 
\quad\triangle = a_{jj}a_{kk}-a_{jk}a_{kj}\tag4.6
$$
\enddemo

We now proceed to reduce the general case to a sum of cases as
described in the preceding remark.

Suppose that  $\pi_1,\pi_2,\dots,\pi_N$, $N=\frac12(n^2-n)$, are
permutation matrices with the two properties (written using the
standard matrix units)
\roster
\item"(4.7a)"   $\pi_{\nu}$  is the matrix of the transposition
$(k,k+1)$ , $k=k_{\nu}$;
\item"(4.7b)"  the product  $\pi_1\pi_2\dots\pi_N$  is the
antidiagonal matrix  $r=\sum e_{j,n+1-j}$.
\endroster

There are various such decompositions of the antidiagonal matrix
$r$.  One such decomposition corresponds to permuting  $(1,2,\dots,n)$
by moving  $1$  to the extreme right in  $n-1$  steps, then moving
$2$  to the position left of  $1$  in  $n-2$  steps, and so on.

\proclaim{(4.9) Definition}  Given  $\pi_{\nu}$  satisfying (4.7), set
$$\align
r_0 &=1,\quad r_{\nu} = \pi_1\pi_2\cdots\pi_{\nu},\quad 1\leq \nu\leq N;\\
U_{\nu} &= r_{\nu}GL^+_0r^{-1}_{\nu};\quad L_{\nu} =
r_{\nu}GL^-r^{-1}_{\nu};\quad D_{\nu} = r_{\nu}B_{\nu}r^{-1}_{\nu},
\endalign
$$
where  $B_{\nu} = \{b\in GL(n); b_{jk}=0$  for  $j\neq k$  unless
$\{j,k\} = \{k_{\nu},k_{\nu} +1)\}$.  Thus the matrices in  $D_{\nu}$  are
block diagonal in the sense that after conjugation by  $r^{-1}_{\nu}$  the
nonzero entries lie in a single  $2\times 2$  block along the
diagonal.
\endproclaim

Note that  $U_{\nu}, L_{\nu},D_{\nu},U_{\nu}+D_{\nu}$,  and  $L_{\nu}
+ D_{\nu}$  are subalgebras of the matrix algebra  $M_n$.  Let
$$
P_{\nu}:M_n \to D_{\nu}\tag4.9
$$
be the projection which commutes with left and right multiplication by
diagonal matrices.  Then
$$
(D_{\nu}+U_{\nu})\cap (D_{\nu}+L_{\nu}) = D_{\nu};\tag4.10
$$
$$
P_{\nu}\quad\text {is an algebra homomorphism on}\quad
D_{\nu}+U_{\nu}\quad\text{and on}\quad D_{\nu}+L_{\nu}. \tag4.11
$$
Note also that because  $r_{\nu-1}$   and  $r_{\nu}$   differ by a single
transposition, we have the identities
$$
D_{\nu} + L_{\nu-1} = D_{\nu} + L_{\nu},\quad D_{\nu} + U_{\nu-1} =
D_{\nu} + U_{\nu}.\tag4.12
$$

\proclaim{(4.13)  Definition}  Suppose the permutation matrices
$\pi_{\nu}$  satisfy (4.7), and suppose  $r^{-1}_{\nu}ar_{\nu}$  is
in  $GL_{\ast}, 0\leq \nu\leq N$.  Then we may factor the
$r^{-1}_{\nu}ar_{\nu}$  as in (4.1) to obtain matrices  $u_{\nu},
\ell_{\nu}$ such that 
$$
au_{\nu} =\ell_{\nu},\quad u_{\nu}\in U_{\nu},\quad \ell_{\nu}\in L_{\nu}.
\tag4.13a
$$
Now define  $v^+_{\nu}, v^-_{\nu}, a^+_{\nu}, a^-_{\nu}$ and $2$-forms
$\Omega_{\nu}$  as follows:
$$\align
v^-_{\nu} &= P_{\nu}u_{\nu-1},\quad v^+_{\nu} = P_{\nu}u_{\nu},\quad
a^-_{\nu} = 
P_{\nu}\ell_{\nu-1},\quad a^+_{\nu} = P_{\nu}\ell_{\nu};\tag4.13b\\
\Omega_{\nu} &=tr[(v^+_{\nu})^{-1}dv^+_{\nu} \wedge (a^+_{\nu})^{-1}
da^+_{\nu}-(v^-_{\nu})^{-1}dv^-_{\nu}\wedge (a^-_{\nu})^{-1}\wedge
da^-_{\nu}]. \tag4.13c
\endalign
$$
\endproclaim

\proclaim{Lemma 4.3}  With  $a$  as in (4.13),
$$
v_- =u_0,\quad a_- =\ell_0,\quad v_+=u_N,\quad a_+=\ell_N\tag4.14.
$$
$$
u_{\nu-1}^{-1}u_{\nu} = \ell_{\nu-1}^{-1}\ell_{\nu} =
(v_{\nu}^-)^{-1}v_{\nu}^+ = (a_{\nu}^-)^{-1}a_{\nu}^+ \tag4.15
$$
In particular,  $a_{\nu}^+(v_{\nu}^+)^{-1}$  and
$a_{\nu}^-(v_{\nu}^-)^{-1}$  have a common value  $a_{\nu}$  in
$D_{\nu}$,  and  $\Omega_{\nu}$  is the pullback at  $a_{\nu}$  of
$\Omega$  under the map  $b\mapsto r^{-1}_{\nu}br_{\nu}$  from
$D_{\nu}$  to  $B_{\nu}\subset GL(n)$.
\endproclaim

\demo{Proof}  That (4.14) holds is clear from the definitions, together
with the assumption (4.7), which gives  $r_N=r$.  The factorizations (4.13a)
imply the first equality in (4.15), since
$\ell_{\nu-1}u^{-1}_{\nu-1} = a=\ell_{\nu}u^{-1}_{\nu}$.
Because of this first identity and (4.12), the common value belongs
to  $(D_{\nu}+L_{\nu})\cap (D_{\nu} + U_{\nu}) = D_{\nu}$.
Therefore we can project and use the property (4.11) to obtain the
remaining identities in (4.15). The final statement is immediate from
conjugation of  $\Omega_{\nu}$  by  $r_{\nu}$.
\enddemo

Note that  $a_{\nu}v_{\nu}^{\pm}=a^{\pm}_{\nu}$, and that  $\Omega_{\nu}$
can be expressed directly in terms of the entries of the  $a_{\nu}$  via
(4.6).  By virtue of the following decomposition theorem, the computation
of  $\Omega$  is reduced to a sum of  $2\times 2$  problems, as in Remark
4.5.  An algorithm for computing  $a_{\nu}$  in terms of  $a$  will be
given below.

\proclaim{Theorem 4.4}  Under the assumption (4.7),  $\Omega$  is the sum
$$
\Omega = \Omega_1+\Omega_2 + \dots + \Omega_N.\tag4.16
$$
\endproclaim

\demo{Proof}  We begin by reversing the reasoning in the proof of
Proposition 3.3 to write  $\Omega$  in the alternative form
$$
\Omega = tr[v_-(dv)v^{-1}_+\wedge a^{-1}da].
$$
From (4.14), (4.15) it follows that
$$
v=v^{-1}_-v_+ = u^{-1}_0u_N=v_1v_2\dots v_N,
$$
where  $v_{\nu}$   is the common value of the matrices in (4.15).  Then
$$
v_-(dv)v_+ = v_-\{\sum(v_1\dots v_{\nu-1})dv_{\nu}(v_{\nu+1}\dots
v_N)\}v^{-1}_+ = \sum u_{\nu-1}(dv_{\nu})u^{-1}_{\nu}.
$$
Note the identity 
$$\align
u_{\nu-1}dv_{\nu}u_{\nu}^{-1} &= u_{\nu}(v_{\nu}^+)^{-1}dv_{\nu}^+
u_{\nu}^{-1}\\
&- u_{\nu-1}(v^-_{\nu})^{-1}dv^-_{\nu}u^{-1}_{\nu-1}.
\endalign
$$
From  $a=\ell_{\nu}u^{-1}_{\nu}$  we find
$$\align
tr\{u_{\nu}(v^+_{\nu})^{-1} dv^+_{\nu} u^{-1}_{\nu} &\wedge a^{-1}da\}\\
= tr\{(v^+_{\nu})^{-1} dv^+_{\nu} &\wedge [(du^{-1}_{\nu})u_{\nu} +
\ell^{-1}_{\nu}d\ell_{\nu}]\}
\endalign
$$
The first term vanishes since  $v^+_{\nu}$  and  $u_{\nu}$  both
belong to  $U_{\nu}$.  For the second term we have
$$\align
tr\{(v^+_{\nu})^{-1} dv^+_{\nu} &\wedge
P_{\nu}(\ell^{-1}_{\nu}d\ell_{\nu})\}\\
= tr\{(v^+_{\nu})^{-1}dv^+_{\nu} &\wedge (a^+_{\nu})^{-1} da^+_{\nu}\}
\endalign
$$
since  $P_{\nu}$  is multiplicative on  $L_{\nu}$.  The second term
in  $tr(u_{\nu-1}dv_{\nu}u^{-1}_{\nu}\wedge a^{-1}da)$  is treated in
the same way, using  $a=\ell_{\nu-1}u^{-1}_{\nu-1}$, and (4.13)
follows immediately.
\enddemo

\demo{Proof of Theorem 4.1}  Choose permutation matrices which lead to
a decomposition (4.16) of  $\Omega$.  According to Lemma 4.2 and
Remark 4.5, each  $\Omega_{\nu}$  has the form  $dp_{\nu}\wedge
dq_{\nu}$, so we obtain the desired representation (4.3) on the dense
set where the   $p_{\nu},q_{\nu}$  are defined.  It follows
immediately that   $\Omega$  is closed and that it has
rank at most  $n^2-n$  everywhere.  The rank is  $n^2-n$  at the identity,
and therefore is  $n^2-n$  on a dense algebraic open set, so the
$dp_{\nu}, dq_{\nu}$   are generically independent.
\enddemo

To compute  $\Omega_{\nu}$, and hence  $p_{\nu},q_{\nu}$, we need to find
$a_{\nu}$.  To obtain  $a_{\nu}$  from  $a$, note that
$$
au^{\#}_{\nu} = \ell^{\#}_{\nu}
\tag4.17
$$
where
$$
u^{\#}_{\nu} =
u_{\nu-1}(v^-_{\nu})^{-1},\qquad  \ell_{\nu}^{\#} =
\ell_{\nu-1}(v_{\nu}^-)^{-1},\tag4.18
$$
and 
$$
P_{\nu}\ell^{\#}_{\nu} = P_{\nu}\ell_{\nu-1}(v^-_{\nu})^{-1} =
a^-_{\nu}(v^-_{\nu})^{-1} = a_{\nu},\quad P_{\nu}u^{\#}_{\nu} = 1 .
$$
The matrix  $r^{-1}_{\nu-1}a_{\nu}r_{\nu-1}$  is block diagonal; its
nontrivial part is the  $2\times 2$  block with entries from
rows and columns  $k_{\nu},k_{\nu +1}$ .
It follows from (4.17), (4.18), and a simple computation that
this block is
$D-CA^{-1}B$,  where  $A$  is  $(k_{\nu}-1)\times (k_{\nu}-1)$, $D$  is 
 $2\times 2$, and
$$
[r^{-1}_{\nu-1}a r_{\nu-1}]_{j,k\leq k_{\nu }+1} = \pmatrix A & B\\C &
D\endpmatrix
$$

The following notational convention will be useful.

\proclaim{(4.19)  Definition}  If  $J$  and  $K$  are two subsets of
$\{1,2,\dots,n\}$  having the same cardinality and  $a$  is in  $M_n$,
then  $m(J;K) = m(J;K;a)$  denotes the determinant of the
corresponding submatrix
$$
m(J;K;a) = \det(a_{jk})_{j\in J,k\in K}.
$$
We set  $m(\emptyset; \emptyset; a)=1$.
\endproclaim
Direct calculation leads to
$$
D-CA^{-1}B = \frac{1}{m(J;J)}\bmatrix m(J,k_{\nu}; J,k_{\nu}) &
m(J,k_{\nu}; J,k_{\nu}+1)\\ m(J,k_{\nu}+1; J,k_{\nu}) & m(J,k_{\nu}+1;
j,k_{\nu}+1)\endbmatrix 
$$ 
where  $J=(1,2,\dots,k_{\nu}-1)$  and  $m(K;L) =
m(K;L;r^{-1}_{\nu-1}ar_{\nu-1})$.  Thus
$$\align
\Omega_{\nu} &= dp_{\nu}\wedge dq_{\nu}\\
p_{\nu} &=
\log\left[\frac{m(J,k_{\nu};J,k_{\nu})m(J,k_{\nu}+1;J,k_{\nu}+1)}
{m(J;J)m(J,k_{\nu},k_{\nu}+1;J,k_{\nu},k_{\nu}+1)}\right]\tag4.20\\
q_{\nu} &=
\log\left[\frac{m(J,k_{\nu}+1;J,k_{\nu})}{m(J,k_{\nu};J,k_{\nu}+1)}\right] 
\endalign
$$

\demo{(4.21)  Example}  For  $n=3$,  we consider the decomposition of  $r$
corresponding to  $(123)\to (312)\to (231)\to (321)$,  i.e.
$$
\bmatrix 0 & 1 & 0\\ 1 & 0 & 0\\ 0 & 0 & 1\endbmatrix \bmatrix 1 & 0
& 0\\ 0 & 0 & 1\\ 0 & 1 & 0\endbmatrix  \bmatrix 0 & 1 & 0\\ 1 & 0 &
0\\ 0 & 0 & 1\endbmatrix  =  \bmatrix 0 & 0 & 1\\ 0 & 1 & 0\\ 1 & 0 &
0\endbmatrix
$$
Let  $A$  be the cofactor matrix  $A = \triangle(a^{-1})^t$, where
$\triangle = \det a$.  Then the corresponding decomposition of   $\Omega$  is
$$\align
\Omega &= d\log (\frac{a_{11}a_{22}}{A_{33}})\wedge d\log
(\frac{a_{21}}{a_{12}})\\ 
&+ d\log (\frac{A_{11}A_{33}}{\triangle a_{22}})\wedge d\log
(\frac{A_{13}}{A_{31}})\\
&+ d\log (\frac{a_{22}a_{33}}{A_{11}})\wedge d\log (\frac{a_{32}}{a_{23}}).
\endalign
$$
\enddemo

We conclude with some symmetry properties of  $\Omega$  which will be
important later.

\proclaim{Proposition 4.7}  The $2$-form  $\Omega$  is odd under the
automorphisms of  $GL(n)$
$$
\Phi_1(a) = (a^{-1})^t,\quad\Phi_2(a) = rar,
$$
(where  $r$  is again the antidiagonal matrix  $\sum e_{j,n+1-j})$,  i.e.
$$
\Phi^{\ast}_j\Omega = -\Omega,\quad j=1,2.
$$
\endproclaim

\demo{Proof}  $\Phi_j(GL^{\pm}) = GL^{\mp}$.  Therefore if
$a_j=\Phi_j(a)$,  the  $\pm$  factors in (4.1) are
$\Phi_j(a_{\mp}),\mathbreak \Phi_j(v_{\mp})$.  It is immediate from
this that $\Phi^{\ast}_2\Omega = - \Omega$.  The result for  $\Phi_1$
makes use also of the identity  $\Phi^{\ast}_1(b^{-1}db) =
-(b^{-1}db)^t$, together with the identity  $tr(\alpha\wedge\beta) =
tr(\alpha^t\wedge\beta^t)$  for matrix-valued $1$-forms.
\enddemo

\head 5.  A symplectic foliation and Poisson bracket on
$GL(n)$; flows \endhead 

We introduce now a foliation of  $GL_{\ast}$   which is naturally
associated to the factorizations (4.1).  As in remark (3.11) we define
diagonal matrices  $\delta_{\pm} = \operatorname{diag}(a_{\pm})$  and
$\delta = \delta^{-1}_-\delta_+$,  and set
$$
a_{\pm}=\delta_{\pm}b_{\pm},\quad b_{\pm}\in GL^{\pm}_0\tag5.1
$$
so that
$$
v^{-1}_-v_+ = a^{-1}_-a_+ = b^{-1}_-(\delta^{-1}_-\delta_+)b_+ =
b^{-1}_-\delta b_+.\tag5.2
$$

In the notation of (4.20) the principal minors of  $a\in GL(n)$  are
$m(J;J;a)$; we abbreviate this to  $m(J;a)$.  If  $a$  is understood, we
may write  $m(J)$.  In particular the upper and lower principal minors
are
$$
d^+_j=d^+_j(a) = m(1,\dots,j);\quad  d^-_j = m(j,\dots,n);\quad  d_0 = 1 =
d_{n+1}.\tag5.3
$$
It follows from (4.1) that  $\delta_+(\delta_-)$  has the same lower
(upper) minors as  $a$, so
$$
(\delta_+)_{jj} = \frac{d^-_j}{d^-_{j+1}};\quad (\delta_-)_{jj} =
\frac{d^+_j}{d^+_{j-1}}, \qquad 1\leq j\leq n.\tag5.4
$$
Therefore  $\delta_+$  and  $\delta_-$  are determined uniquely by
$\delta = \delta^{-1}_-\delta_+$,  together with the quotients
$$
\varphi_j(a) = d^+_j(a)/d^-_{j+1}(a),\qquad 1\leq j\leq n.\tag5.5
$$
Note that  $\varphi_n(a) = \det a$,  and that the  $\varphi_j$
determine the product $\delta_+\delta_-$.

\proclaim{Theorem 5.1}  The foliation of  $GL_{\ast}$  by the
functions  $\varphi_j$  is symplectic for  $\Omega$.  Each
leaf  $\{a:\varphi_j(a) = c_j, \quad 1\leq j\leq n\}$
is parametrized by  $V_{\ast} = \{(v_+,v_-)\in GL^-_0\times
GL^+_0:v^{-1}_-v_+\in GL\}$;  the pullback of  $\Omega$  to a leaf is
generically of rank  $n^2-n$; the pullback from a leaf to  $V_{\ast}$
is independent of the choice of leaf and is given by
$$
tr[v^{-1}_+dv_+\wedge \{b^{-1}_+db_+ +
b^{-1}_+(\delta^{-1}d\delta)b\}
-v^{-1}_-dv_-\wedge \{b^{-1}_-db_- -
b^{-1}_-(\delta^{-1}d\delta)b_-)\}]\tag5.6
$$
where  $b_{\pm}$  and  $\delta$  are determined from  $(v_+,v_-)$  by
the factorization 
$$
v^{-1}_-v_+ = b^{-1}_-\delta b_+, \quad b_{\pm}\in GL^{\pm}_0, \quad
\delta = \operatorname{diag}(\delta).\tag5.7
$$
\endproclaim

\demo{Proof}  Starting with  $a$  in  $GL_{\ast}$, we define
$b_{\pm}, \delta_{\pm}$,  as above and again set  $\delta =
\delta^{-1}_-\delta_+$.  Let  $\eta = \delta_-\delta_+$.  Then
$$
\delta^{-1}_{\pm}d\delta_{\pm} = \eta^{-1}d\eta\pm
\delta^{-1}d\delta;\quad a^{-1}_{\pm}da_{\pm} = b^{-1}_{\pm}db_{\pm}
+ b^{-1}_{\pm}(\delta^{-1}_{\pm}d\delta_{\pm})b_{\pm}.
$$
Consequently  $\Omega$  is given by the sum of (5.6) and
$$
tr[v^{-1}_+dv_+\wedge b^{-1}_+(\eta^{-1}d\eta)b_+ -
v^{-1}_-d\nu_-\wedge b^{-1}_-(\eta^{-1}d\eta)b_-].\tag5.8
$$
The pullback of  $d\eta$  to the leaves of the foliation determined by the
functions  $\varphi_j$  vanishes, since these functions determine
$\eta$, so the pullback of  $\Omega$  to the leaf is given by (5.6).
Now  $b_{\pm}$  and  $\delta$  are determined from  $(v_+,v_-)$  in
$V_{\ast}$  by (5.2), so the pullback of  $\Omega$  to
$V_{\ast}$  is leaf-independent and given also by (5.6).  These
pullbacks have rank  $\leq\dim(V_{\ast}) = n^2-n$
everywhere and rank  $n^2-n$  at the unique diagonal element in
a given leaf, so they have rank  $n^2-n$  generically.
\enddemo

The symplectic foliation gives rise to a Poisson structure on
$GL_{\ast}$.  In fact there is a Poisson bracket  $(~,~)_L$   on each
leaf  $L$  corresponding to the pullback of  $\Omega$  to  $L$, and
this may be extended to a (degenerate) Poisson bracket for functions
on  $GL_{\ast}$,  characterized by
$$
(f,g)\mid_L = (f\mid_L,g\mid_L)_L; \quad (f,\varphi_j) = 0, \quad
1\leq j\leq n.\tag5.9
$$
Equivalently, if  $\Omega = \Sigma dp_{\nu}\wedge dq_{\nu}$   as in
Theorem 4.1, then
$$\multline
(p_{\mu},q_{\nu}) = \delta_{\mu\nu},\quad (p_{\mu},p_{\nu}) = 0 =
(q_{\mu},q_{\nu}), \quad 1\leq \mu,\nu \leq (n^2-n);\\
(f,\varphi_j) = 0, \quad 1\leq j\leq n. \endmultline\tag5.10
$$
Functions such as the  $\varphi_j$  which Poisson commute with all
functions are sometimes called  \it Casimirs\rm.

This Poisson bracket was computed for the standard coordinates of
$GL(n)$  by Lu \cite{Lu} in the cases  $n=2$, $n=3$; Lu also
conjectured the general form below and pointed out the connection
with the classical limit of a quantum version due to Drinfeld
\cite{Dr}, as noted in section 1.

\proclaim{Proposition 5.2}  The Poisson bracket (5.9) is odd under
the automorphisms  $\Phi_1,\Phi_2$  of Proposition 4.7, i.e.
$$
(f\circ\Phi_j,g\circ\Phi_j) = -(f,g)\circ\Phi_j,\quad j=1,2.
$$
\endproclaim

\demo{Proof}  The Casimirs  $\varphi_k$  satisfy
$\varphi_k\circ\Phi_j=\varphi_{n+1-k}$,  so  $\Phi_j$  maps leaves to
leaves.  Proposition 4.7 implies that the pullback under  $\Phi_j$
to a leaf  $L$  of  $\Omega$  on  $\Phi_j(L)$  is  $-\Omega$.
Therefore the pushforward of the Poisson bracket  $(~~,~~)_L$  is
$-(~~,~~)_{\Phi_{i}(L)}$,  and (5.9) implies the desired result.
\enddemo

\proclaim{Theorem 5.3}  The Poisson bracket on  $GL_{\ast}$  induced
by the $2$-form  $\Omega$  and the foliation by functions
$(\varphi_j)$  extends to the full matrix space  $M_n$  and is given by
the following bracket relations between the coordinate functions
$a_{jk}$:
$$
(a_{jk},a_{\ell m}) = \frac14 [sgn(\ell -j) - sgn(m-k)]a_{jm}a_{\ell
k},\tag5.11
$$
where again  $sgn(0) = 0$.
\endproclaim 

\demo{Proof}  We show first that the calculation can be reduced to the
cases  $n<4$.  For a fixed  $\kappa$, $1\leq\kappa\leq n$,  consider
the map  $M_n\to M_{n-1}$  obtained by omitting the  $\kappa$-th row and
column.  We claim that the pushforward to  $M_{n-1}$  of the Poisson
structure on   $M_n$ coincides with the structure on  $M_{n-1}$
itself, i.e. the Poisson bracket of coordinate functions  $a_{jk},
a_{\ell m}$  on  $M_n,\quad i,j,k,\ell\neq\kappa$,  is the same
as that obtained by considering them as functions on  $M_{n-1}$.  To verify
this claim we take the decomposition (4.7) of the antidiagonal matrix
$r\in M_n$  obtained from the following three sets of permutations: the
first set takes  $(1,2,\dots,n)$  to
$(1,\dots,\kappa-1,\kappa+1,\dots,n)$  in  $n$-steps; the second set
takes us to  $(n,n-1,\dots,1,\kappa)$  in  $\frac12(n-1)(n-2)$
steps; the third set takes us to  $(n,n-1,\dots,1)$  in  $(\kappa-1)$
steps.  The corresponding additive decomposition of $\Omega$  then
takes the form 
$$
\Omega = \Omega' + \Omega'' + \Omega''' \tag5.12
$$
with the obvious notational convention.  According to the prescription
after Lemma 4.5,  $\Omega'' = \Sigma dp_{\nu}\wedge dq_{\nu}$, where
the  $p_{\nu}, q_{\nu}$  are functions of the matrix  $b$  in
$M_{n-1}$  which corresponds to  $a$  in  $M_n$  by the map above;
moreover   $\Omega''$  has exactly the same form as  $\Omega_{n-1}$.
To complete the verification of our claim it is, therefore, sufficient
to show that the foliation functions for  $b$  Poisson commute with all
entries of  $b$  when these are considered as functions of  $a$.  We know
that the  $p_{\nu}$   which correspond to  $\Omega'$  and to
$\Omega''$  in the decomposition (5.12) commute with all  $p_{\nu},
q_{\nu}$  from  $\Omega''$  and with each other, so it is enough to
show that the foliation functions of  $b$ in  $M_{n-1}$   are
computable from the  $p_{\nu}$  in  $\Omega'$  and  $\Omega'''$, together
with the foliation functions on  $M_n$.  The  $p_{\nu}$
corresponding to  $\Omega'$  are
$$
\log(g_{k+1}/g_k), \quad\kappa\leq k < n;\quad g_k =
m(1,\dots,\hat{\kappa},\dots,k)/m(1,\dots,k)
$$
where again  $m(\dotso)$  denotes the principal minor of  $a$  based
on the indicated rows and columns.  Similarly, the  $p_{\nu}$
associated to  $\Omega'''$  are
$$
\log(h_{k+1}/h_k),\quad 1\leq k < \kappa;\quad
h_k=m(k,\dots,n)/m(k,\dots,\hat{\kappa},\dots,n).
$$
Modulo the foliation functions on  $M_n$,  the  $g_k$  and  $h_k$
can be determined from the  $p_{\nu}$.  Again modulo the
foliation functions on  $M_n$, the  $g_k$'s  and $h_k$'s  are
equivalent to the set of functions
$$\align
m&(1,\dots,\hat{\kappa},\dots,k)/m(k+1,\dots,n), \quad\kappa\leq k\leq n;\\
m&(1,\dots,k-1)/m(k,\dots,\hat{\kappa},\dots,n), \quad 1\leq k <\kappa.
\endalign
$$
These are precisely the foliation functions of  $b$  as an element of
$M_{n-1}$.  This completes the proof that the two Poisson structures
coincide on  $b$.

Suppose now that  $a_{ij}, a_{k\ell}$  are any two coordinate functions on
$M_n$.  Repeated use of the preceding argument shows that their Poisson
bracket can be computed by taking them to be functions on  $M_p$,  with
$p$  the cardinality of  $\{i,j,k,\ell\}$.  Thus the computation is reduced to
the cases  $M_1$  (trivial),  $M_2,M_3,M_4$.  The complete computation is
tedious, so we merely indicate a few representative cases.  Recall
that a Poisson bracket is a derivation for each of its arguments.

For  $n=2$  the foliation functions are  $a_{11}/a_{22},\quad \Delta
= \det a$.  From this fact and (4.4) we deduce
$$\align
(a_{11},a_{22}) &= (a_{11}, a_{11}(a_{22}/a_{11})) = (a_{11},
a_{11})a_{22}/a_{11} = 0;\\
(a_{11}, a_{12}a_{21}) &= (a_{11}, a_{11}a_{22}-\Delta) = 0.
\endalign
$$
Also,  $p=\log(a_{11}a_{22}/\Delta)$, $q=\log(a_{21}/a_{12})$  and
$(p,q) = 1$,  so
$$\align
(a_{11}a_{22}, a_{21}/a_{12}) &= a_{11}a_{22}a_{21}/a_{12};\\
(a_{11}, a_{21}/a_{12}) &= (a^2_{11}, a_{21}/a_{12})/2a_{11} =
(a_{11}a_{22}a_{11}/a_{22}, a_{21}/a_{12})/2a_{11}\\
&= (a_{11}a_{22}, a_{21}/a_{12})/2a_{22} = a_{11}a_{21}/2a_{12}.
\endalign
$$
Therefore
$$
(a_{11}, a_{21}) = \frac{1}{2a_{21}}(a_{11},a^2_{21} = a_{12}(a_{11},
a_{21}/a_{12})/2 = \frac14 a_{11}a_{21}
$$
Because of the symmetries in Proposition 4.7,  $(a_{22},a_{21}) = -
(a_{11},a_{21})$  and so on.  Also
$$\align
(a_{21},a_{12}) &= (a_{21},a_{12}a_{21})/a_{21} = (a_{21},
-a_{11}a_{22})/a_{21}\\
&= (a_{11}a_{22}a_{11}/a_{22}, a_{21})a_{21}/a_{11}a_{12} = 2(a_{11},
a_{21})a_{22}/a_{21} = a_{11}a_{22}/2.
\endalign
$$

For  $n=3$, reduction to  $n=2$  gives all brackets such as
$(a_{22},a_{12}), (a_{23}, a_{33})$.  Let  $A=(\det a)(a^{-1})^t$  be
the cofactor matrix.  The foliation functions are
$$
\det a,\quad a_{11}/A_{11}, \quad a_{33}/A_{33}.
$$
The decomposition (4.21) implies, therefore, that
$$\multline
0 = (a^2_{22}a_{11}a_{13}/A_{11}A_{33}, A_{13}/A_{31}) =
2a_{22}a_{11}a_{33}(a_{22}, A_{13}/A_{31})/A_{11}A_{33}\\
(a_{22}, A_{13}A_{31}) = (a_{22}, A_{11}A_{33}-a_{22}\det a) =
A_{11}A_{33}(a_{22}, a_{11}a_{33})/a_{11}a_{33} = 0.\endmultline
$$
Therefore  $(a_{22}, A_{12}) = 0 = (a_{22}, A_{31})$,  and the
$2\times 2$  results allow one to calculate  $(a_{22}, a_{31})$  and
$(a_{22}, a_{13})$  from these identities.  Similar computations
yield all the  $3\times 3$  brackets, though for some it is
convenient to replace the decomposition in (4.21) with the 
decomposition obtained from factoring  $r$  by means of  $(123)\to
(132)\to (312) \to (321)$.

The case  $n=4$  is similar, and again brackets like $(a_{12},
a_{24})$  are known from the  $3\times 3$  computation.  This
completes our sketch of the proof of Theorem 5.3.

We consider now the  $(n-1)$-parameter family of flows in  $M_n$:
$$
a(t) = \exp(t\mu)a(0)\exp(-t\mu),\quad \mu\quad\text{diagonal},\quad
tr\mu = 0.\tag5.13
$$
This conjugation preserves the principal minors of  $a$, so the flow
preserves  $GL(n)$, $GL_{\ast}$, and the leaves of the foliation.  The
factorization (4.1) is also preserved, so the flow preserves the
$2$-form  $\Omega$.  Therefore these flows are Hamiltonian. 
\enddemo

\proclaim{Theorem 5.4}  The Hamiltonian function for the flow (5.13) is
$tr[\mu\log\delta], \quad \delta = \delta^{-1}_-\delta_+$,  and it is
a linear combination of the functions  $p_{\nu}$  of (4.20).
\endproclaim

\demo{Proof}  We use the additive decomposition of Theorem 4.4, with
the factorization of  $r$  described after (4.7).  If  $\pi_{\nu}$ is
associated with the interchange of  $j$  and  $k$, then (4.20)
implies that under the flow (5.13),  $\dot q_{\nu} = \mu_k-\mu_j$.
Therefore the Hamiltonian for (5.13) is a linear combination of the
$p_{\nu}$.  The  $p_{\nu}$  themselves are logarithms of quotients of
principal minors of  $a$:
$$
m_{jk} = m(j,j+1,\dots,k),\quad j\leq k
$$
The term  $m_{jk}$  occurs in the \it numerator \rm of
$\exp(p_{\nu})$  when  $\pi_{\nu}$  is associated to the interchange
of  $j-1$  with  $k$  or of  $j$  with  $k+1$;  the term  $m_{jk}$
occurs in the \it denominator \rm when   $\pi_{\nu}$  is associated to 
the interchange of  $j-1$  with  $k+1$  or of  $j$  with  $k$.  Thus
the total weight attached to  $\log m_{jk}$  in the Hamiltonian for
the flow (5.13) is
$$
(\mu_k-\mu_{j-1}) + (\mu_{k+1}-\mu_j) - (\mu_{k+1}-\mu_{j-1}) -
(\mu_k-\mu_j),
$$
with the convention that if  $j-1=0$  or  $k+1=n$, the corresponding
term in parentheses is omitted.  Thus  $\log m_{jk}$  has weight zero
unless  $j=1$  or  $k=n$,  and the Hamiltonian for (5.13) is
$$\multline
\sum^{n-1}_{k=1}(\mu_{k+1}-\mu_k)\log m_{1k} +
\sum^{n-1}_{j=1}(\mu_{j+1}-\mu_j)\log m_{jn}\\
= \sum^n_{j=1}\mu_j\log [d^-_j(a)d^+_{j-1}(a)/d^-_{j+1}(a)d^+_j(a)] =
\sum^n_{j=1}\mu_j\log\delta_{jj}.
\endmultline
$$
\enddemo

\demo{(5.14)  Remark}  Theorem 5.2 and its proof show that the flows
(5.13) are completely integrable in the classical sense: they are an
$(n-1)$-parameter family of commuting Hamiltonian flows in each  $n^2-n$
dimensional symplectic leaf of the foliation, and are part of the
$\frac12(n^2-n)$-parameter family of commuting flows generated by the
$p_{\nu}$.  Note that the flows (5.13) are the only members of the
larger family which are linear as flows on the full matrix algebra
$M_n$; in fact the generator of a linear flow which commutes with all
the flows (5.13) must have each matrix unit  $e_{jk}, j\neq k$, as an
eigenvector and if such a flow leaves all the  $p_{\nu}$  invariant
it can be shown to be included among the flows (5.13).
\enddemo

\head 6.  Proofs of Theorems 3.7 and 3.8: Darboux coordinates for
scattering data \endhead

Up to a trivial normalization, the functions  $p_{\nu}, q_{\nu}$  of
Theorem 4.1  are the functions of Theorem 3.7 and Theorem 3.8.   To
see this, we return to the notation introduced before Theorem 3.7.
Observe that if  $f,g$  belong to  $C(SL(n))$  then in view of
Proposition 3.4 the distribution-valued Poisson bracket can be
decomposed as 
$$
\{f_{\xi}, g_{\eta}\}_S = [f,g](s(\xi))\delta(\xi-\eta) + \langle
f,g\rangle (s(\xi),s(\eta))  p.v. \frac{1}{\xi-\eta}\tag6.1
$$
where  $[~,~]$  and  $\langle~,~\rangle$  have the following properties.
\roster
\item"(6.2)"  The map  $f,g\mapsto [f,g]$ is an alternating bilinear
map from  $C^{\infty}(SL(n))\times C^{\infty}(SL(n))$  to
$C^{\infty}(SL(n))$  which is a derivation in each variable:  
$$\align
[fg,h] &= f[g,h] + g[f,h]\\
[f,gh] &= g[f,h] + h[f,g].
\endalign
$$
\item"(6.3)"  $[a_{jk},a_{\ell m}] = \pi i a_{jm}a_{\ell
k}[sgn(\ell-j)-sgn(m-k)]$. 
\item"(6.4)"  The map  $f,g\to\langle f,g\rangle$  is a symmetric
bilinear map from  $C^{\infty}(SL(n))\times C^{\infty}(SL(n))$  to
$C^{\infty}(SL(n)\times SL(n))$  such that
$$\align
\langle fg,h\rangle(s,s') &= f(s)\langle g,h\rangle(s,s') +
g(s)\langle f,g\rangle(s,s')\\
\langle f,gh\rangle(s,s') &= g(s')\langle f,h\rangle(s,s') +
h(s')\langle f,g\rangle(s,s'). 
\endalign
$$
\item"(6.5)"  $\langle a_{jk},a_{\ell m}\rangle(s,s') =
a_{jk}(s)a_{\ell m}(s')[\delta_{j\ell} - \delta_{km}]$.
\endroster
Here the  $a_{jk}$  are the coordinate functions and  $s,s'$  are points of
$SL(n)$.

The properties (6.2), (6.3) imply that the bracket  $[~,~]$  is
precisely  $4\pi i(~,~)$,  where  $(~,~)$  is the Poisson bracket
(5.11);  this corresponds to the fact that  $\Omega_S$  is the direct
integral of $4\pi i \Omega$.  Consequently we may replace the
functions  $p_{\nu}$  and  $q_{\nu}$  of section 4  by the
renormalized versions 
$$
\frac12 p_{\nu},\quad \frac{1}{2\pi i}q_{\nu}\tag6.6
$$
to obtain
$$
\{p_{\mu,\xi},q_{\nu,\eta}\}_S = \delta_{\mu\nu}\delta(\xi-\eta) +
\langle p_{\mu}, q_{\nu}\rangle (s(\xi),s(\eta))  p.v.
\frac{1}{\xi-\eta}\tag6.7
$$
To complete the proof of Theorem 3.7  we must show
$$\align
\langle p_{\nu},f\rangle &= 0,\quad \text{all}\quad
f\quad\text{in}\quad C^{\infty}(SL(n)),\quad\text{all}\quad
\nu;\tag6.8\\
\langle q_{\mu},q_{\nu}\rangle &= 0,\quad\text{all}\quad\mu,\nu.\tag6.9
\endalign
$$

\proclaim{(6.10)  Definition}  Given subsets  $J,J',K,K'$  of
$\{1,2,\dots,n\}$,  set
$$
\varepsilon(J,K;J',K') = \text{card}(J\cap J') - \text{card}(K\cap
K').
$$
\endproclaim

\proclaim{Lemma 6.1}  The bracket  $\langle~,~\rangle$   between
minors of  $s$  satisfies 
$$
\langle m(J;K),m(J';K')\rangle(s,s') =
\varepsilon(J,K;J',K')m(J;K;s)m(J';K';s') \tag6.11
$$
\endproclaim

\demo{Proof}  The case when  $J,K,J',K'$  all have cardinality  $1$ is
immediate from 6.5.  The general case follows by expanding the
determinants and using the derivation property 6.4.
\enddemo

We can now prove (6.8) and (6.9), using (4.20).  Each  $p_{\nu}$  is the
logarithm of a term  $m(J)m(J')/m(K)m(K')$,  where each element  $j$  of
$(1,2,...,n)$  occurs with the same frequency in the pair of sets $J,J'$
as in the pair  $K,K'$.  From the derivation property (6.4) we deduce
that (6.11) is equivalent to
$$
\langle\log m(J;K),\log m(J';K')\rangle =
\varepsilon(J,K;J',K').\tag6.12
$$
Each coordinate function  $a_{jk}$  is itself a minor, and therefore
(6.12) implies 
$$\align
\langle\log &\frac{m(J)m(J')}{m(K)m(K')},\log a_{jk}\rangle\\
&=\varepsilon(J,J;j,k) + \varepsilon(J',J';j,k) -
\varepsilon(K,K;j,k) -\varepsilon(K',K';j,k)\\
&= 0.
\endalign
$$

This proves (6.8).  To prove (6.9) we note that according to (4.20),
each  $q_{\nu}$  has the form  $\log m(J;K)/m(K;J)$, which we
abbreviate slightly as  $\log m(JK)/m(KJ)$.  Again we deduce from
(6.12) that

$$\align
\langle\log &\frac{m(JK)}{m(KJ)},
\log\frac{m(J'K')}{m(K'J')}\rangle\\
&= \varepsilon(J,K;J',K') - \varepsilon(J,K;K',J) -
\varepsilon(K,J;J',K') + \varepsilon(K,J;K',J')\\
&= 0
\endalign
$$
because  $\varepsilon(J,K;J'K') = - \varepsilon(K,J;K',J')$,  and
$\varepsilon(J,K;K',J') = - \varepsilon(K,J;J',K')$.  This proves 
(6.9).  For the injectivity property (3.21) we need the foliation
functions  $\varphi_j$  in addition to  $p_{\nu}, q_{\nu}$.  By
Theorem 5.4, entries of  $\delta^{-1}_-\delta_+$  are linear
combinations of the  $p_{\nu}$.  As in Remark 3.11, $\delta_-$ and
$\delta_+$  are Riemann-Hilbert factors of  $\delta^{-1}_-\delta_+$.
Finally, (5.4) and (5.5) determine the $\varphi_j$  from
$\delta_+,\delta_-$. \qed

\demo{Proof of Theorem 3.8}  It is a well-known fact that the
Hamiltonian for the flow (2.8) is the negative of the coefficient of
$z^{-k-1}$  in the asymptotic expansion of  $tr\mu\log\delta(z)$  as
$z\to\infty$,  where 
$$
\delta(z) = \lim\limits_{x\to +\infty}\psi(x,z)\exp(-xzJ), \quad z\in\bold
C\backslash\bold R.
$$
See \cite{S}, \cite{BC3}.  Now  $\delta$  is piecewise holomorphic
with limits  $\delta_{\pm}$  on  $\bold R$,  so the Hamiltonian can
be expressed in the form
$$\align
\frac{1}{2\pi i} &\int_{\bold R}\xi^k[tr\mu\log\delta_+(\xi) -
tr\mu\log\delta_-(\xi)]d\xi\\
&= \frac{1}{2\pi i}\int_{\bold
R}\xi^ktr\mu\log[\delta^{-1}_-(\xi)\delta_+(\xi)]d\xi.
\endalign
$$
According to Theorem 5.4, this integral is a linear combination of
the integrals (3.25).
\enddemo

\head 7.  Coordinates and flows on  $SU(3)$ \endhead

Formula (4.2) defines a complex $2$-form on the intersection of
$GL = GL(n,\bold C)$  with the real submanifold   $GL_{\ast}\cap
U(n)$.  The automorphism  $a\mapsto (a^{-1})^{\ast}$  takes
$GL^{\pm}$  to  $GL^{\mp}$.  Therefore, in the factorizations
(4.1), we have
$$
a_{\pm} = (a^{-1}_{\mp})^{\ast},\quad v_{\pm} =
(v^{-1}_{\mp})^{\ast},\quad a\in U(n).\tag7.1
$$
As in the proof of Proposition 4.7 we can conclude that  $\Omega =
-\overline{\Omega}$  on  $GL_{\ast}\cap U(n)$.

\proclaim{Theorem 7.1} (a) $i\Omega$  is a closed real $2$-form on
$GL_{\ast}\cap U(n)$. 

     (b) The foliation of  $GL_{\ast}\cap U(n)$  induced by the
foliation of  $GL_{\ast}$  in section 5 has leaves with real
dimension  $n^2-n$.  The $2$-form generically has rank  $n^2-n$  on
each leaf, so the foliation is again symplectic.

     (c) If  $\mu$  is a real diagonal matrix with  $tr(\mu) = 0$,
the flow 
$$
a(t) = e^{it\mu}a(0)e^{-it\mu}\tag7.2
$$
is Hamiltonian in  $GL_{\ast}\cap U(n)$  with real Hamiltonian function
$tr(\mu\log\delta)$,  where again  $\delta = \delta^{-1}_-\delta_+,
\quad\delta_{\pm} = \operatorname{diag}(a_{\pm})$.
\endproclaim

\demo{Proof}  Part (a) follows from the preceding remarks.  For part
(b), note that (7.1) implies
$$
\delta_{\pm}(a) = [\delta_{\mp}(a)^{-1}]^{\ast}, \quad a\in U(n).\tag7.3
$$
Because of (5.3) and (5.4), (7.3) implies that the foliation functions
(5.5) take values in  $\{\vert z\vert = 1\}$, so the induced
foliation is defined by the  $n$ independent real functions
$\operatorname{arg}\varphi_j$.  The $2$-form  $i\Omega$  has
(real) rank  $n^2-n$  at the unique diagonal element in each leaf, so it
has rank  $n^2-n$  generically.  Finally, the Poisson bracket associated
to this symplectic foliation is  $-i(~,~)$,  where  $(~,~)$  is the
Poisson bracket of section 5 restricted to  $U(n)$.  According to Theorem
5.3, therefore, the Hamiltonian for the flow (7.2) is  $tr(\mu\log\delta)$,
and according to (7.3) $\delta = \delta^{-1}_-\delta_+ =
\delta^{\ast}_+\delta_+$  is real. \enddemo

The Darboux coordinates  $p_{\nu},q_{\nu}$  constructed in section 4  are
not real when specialized to  $U(n)$  (and suitably normalized) except
for  $n=2$.  We show in this section that real Darboux coordinates can
be chosen in  $SL(3,\bold C)$  in such a way that:  (a)  the
restrictions to  $SU(3)$  are real;  (b)  linear combinations of the
$p_{\nu}$  still include the Hamiltonians for the flows (7.1).  As we
shall see, the third Hamiltonian flow which commutes with the
$2$-parameter family of linear flows (7.2) is not linear on  $M_3$.

We begin our discussion by recalling the Darboux coordinates for
$SL(3)$  in example (4.21).  Corresponding to the new normalization
$i\Omega$,  we take these to be
$$\align
\tilde p_1 &= \log(a_{11}a_{22}/A_{33}),\quad \tilde p_2 =
\log(A_{11}A_{33}/a_{22}), \quad \tilde p_3 = \log(a_{22}a_{33}/A_{11});\\ 
\tilde q_1 &= i\log (a_{21}/a_{12}), \quad \tilde q_2 =
i\log(A_{13}/A_{31}), \quad \tilde q_3 = i\log(a_{32}/a_{23}).
\endalign
$$
Here again the  $A_{jk}$  are the entries of the cofactor matrix
$(\det a)(a^{-1})^t$.  Therefore
$$
A_{jk} = \bar a_{jk}\quad\text{for}\quad a\quad\text{in}\quad
SU(n).\tag7.5
$$
It is convenient to make a preliminary linear canonical transformation
to new Darboux coordinates
$$\align
p_1 &=\tilde p_1 + \tilde p_2 = \log(a_{11}A_{11})\tag7.6\\
p_2 &= \tilde p_1 + \tilde p_2 + \tilde p_3 =
\log(a_{11}a_{22}a_{33})\\
p_3 &= \tilde p_2 + \tilde p_3 = \log(a_{33}A_{33})\\
q_1 &= \tilde q_2 - \tilde q_3 = i\log(a_{23}A_{13}/a_{32}A_{31})\\
q_2 &= \tilde q_1 - \tilde q_2 + \tilde q_3 = i\log
(a_{21}A_{31}a_{32}/a_{12}A_{13}a_{23})\\ 
q_3 &= -\tilde q_1 + \tilde q_2 = i\log(A_{13}a_{12}/A_{31}a_{21}).
\endalign
$$
Then  $p_1$  and  $p_3$  are real on  $SU(3)$, but  $p_2$ is not.  It
follows either from direct calculation using the Poisson bracket (5.11) or
from the first symmetry in Proposition 5.2 that
$\log(A_{11}A_{22}A_{33})$ is also in involution with  $p_1$  and
$p_2$,  so we may take
$$
I_j=a_{jj}A_{jj},\quad j=1,2,3,\tag7.7
$$
as the action variables for a new set of Darboux coordinates.  The
corresponding angle variables  $\Theta_1,\Theta_2,\Theta_3$  are then
obtained by the classical Liouville method; cf \cite{W}.  We briefly
recall the method:  in principle one solves the equations
$$
p_j=f_j(I,q) = f_j(I_1,I_2,I_3,q_1,q_2,q_3),\quad j=1,2,3.\tag7.8
$$
Because the  $I_j$  are independent and in involution, it follows
that  $\partial f_j/\partial p_k = \partial f_k/\partial q_j$  on the
level surfaces  $\{I_j=c_j,~~j=1,2,3\}$.  Therefore there is a generating
function  $S(q,I)$  such that
$$
\partial S/\partial q_j=f_j(I,q) = p_j,\quad j=1,2,3.\tag7.9
$$
Putting  $\Theta_j=\partial S/\partial I_j$,  it follows that
$$
dS = \Sigma p_jdq_j + \Sigma\Theta_jdI_j
$$
and therefore
$$
i\Omega = \Sigma dp_j\wedge dq_j = \Sigma dI_j\wedge
d\Theta_j.\tag7.10
$$
By (7.5) the functions  $I_j$  are real on  $SU(3)$  and we may
replace the  $\Theta_j$  in (7.10) by their real parts (if necessary)
to obtain real Darboux coordinates.  In the remainder of this section
we show that the angle variables  $\Theta_j$  are elliptic functions.

We turn to equation (7.8).  We already have  $p_j=\log I_j,\quad
j=1,3$.  The remaining equation can be obtained, in principle, once
we have a nontrivial identity involving  $p_2, I_j, q_j$.

\proclaim{Proposition 7.2}  Let  $\zeta =
e^{p_{2}}=a_{11}a_{22}a_{33}, ~~I_0=I_2-I_1-I_3+1$.  The following
identity holds on  $SL(3,\bold C)$:
$$
(\zeta I_0 + 2I_1I_3)^2 +
4\cos^2(\frac12 q_2)(\zeta-I_1)(\zeta-I_3)(\zeta-I_1I_3) = 0.\tag7.11
$$
\endproclaim

\demo{Proof}  We make extensive use of the identities defining the
cofactors:  $A_{12}=a_{23}a_{31}-a_{21}a_{33}$  and so on, and of the
corresponding identities coming from the inverse matrix
$A:a_{12}=A_{23}A_{31}-A_{21}A_{33}$  and so on.  In the following
computation each term in braces is replaced by its expression
obtained from such identities in order to pass to the next in the
sequence of identities.
$$\align
I_2 &=a_{22}A_{22}=a_{11}a_{22}a_{33} -
\{a_{22}a_{13}\}\{a_{22}a_{31}\}/a_{22}\\
&= a_{11}a_{22}a_{33} - [\{A_{13}A_{31}\} +
\{a_{12}a_{21}\}\{a_{23}a_{32}\} - (a_{12}A_{13}a_{23} +
a_{21}A_{31}a_{32})]/a_{22} \\
&= I_1+I_3+1-2A_{11}A_{13}/a_{22} + (a_{12}A_{13}a_{23} +
a_{21}A_{31}a_{32})/a_{22} 
\endalign
$$
or
$$\align
a_{12}A_{13}a_{23} + a_{21}A_{31}a_{32} &= a_{22}I_0 +
2A_{11}A_{33}\tag7.12\\
&= (\zeta I_0 + 2I_1I_3)/a_{11}a_{33}.
\endalign
$$
On the other hand
$$\align
(a_{12}A_{13}a_{23})(a_{21}A_{31}a_{32}) &=
(a_{12}a_{21})(A_{13}A_{31})(a_{23}a_{32})\tag7.13\\
&= (a_{11}a_{22}-A_{33})(A_{11}A_{33}-a_{22})(a_{22}a_{33}-A_{11})\\
&= \frac{1}{a_{33}}(\zeta -
I_3)\frac{1}{a_{11}a_{33}}(I_1I_3-\zeta)\frac{1}{a_{11}}(\zeta-I_1) 
\endalign
$$
Let  $f=a_{12}A_{13}a_{23}$  and  $g=a_{21}A_{31}a_{32}$.  Then
$$\align
(f+g)^2 &= fg[\sqrt{f/g} + \sqrt{g/f}]^2 = fg[\exp(-\frac12 iq_2) +
\exp(\frac12 iq_2)]^2\tag7.14\\ 
&= 4fg\cos^2(\frac12 q_2).
\endalign
$$
Combining (7.12), (7.13), and (7.14), we obtain (7.11).

Returning to the generating function  $S$,  we have
$$
\partial S/\partial q_1 = p_1 = \log I_1,\quad\partial S/\partial q_3 =
p_3 = \log I_3,\quad \partial S/\partial q_2 = p_2=\log\zeta
$$
so
$$
S(I,q) = q_1\log I_1 + q_3\log I_3 +
\int^{q_{2}(\zeta,I)}\log\zeta(I,u)du.
$$
The angle variables corresponding to the action variables  $I_j$  are
$$\align
\Theta_j &= \frac{\partial S}{\partial I_j} = \frac{q_j}{I_j} +
\int^{q_{2}}\frac{1}{\zeta}\frac{\partial\zeta}{\partial I_j}du,
j=1,3;\\
\frac{\partial S}{\partial I_2} &=
\int^{q_{2}}\frac{1}{\zeta}\frac{\partial\zeta}{\partial I_2}du. 
\endalign
$$
We rewrite (7.11) in the form
$$
\Phi(\zeta,I,u) = F(\zeta,I) + \cos^2(u/2)G(\zeta,I) = 0\tag7.16
$$
to define  $\zeta$  or  $u$  in terms of  $u$  or  $\zeta$  and  $I$.
Thus
$$
\frac{\partial\zeta}{\partial I_j} = -\frac{\partial\Phi/\partial
I_j}{\partial\Phi/\partial\zeta} =
\frac{(\partial\Phi/\partial I_j)(\partial\Phi/\partial
u)}{\partial\zeta/\partial u}.  
$$
In particular, on the surface  $\Phi=0$,
$$\align
\frac{\partial\Phi}{\partial I_2} &= \frac{\partial F}{\partial I_2} =
2\zeta\sqrt{F};\quad \frac{\partial\Phi}{\partial u} =
-\sin(\frac12 u)\cos(\frac12 u)G;\\
\cos(\frac12 u) &= i[F/G]^{\frac12};\quad \sin(\frac12 u) = [1+F/G]^{\frac12};
\endalign
$$
so
$$
\Theta_2 = \frac{\partial S}{\partial I_2} =
i\int^{q_{2}}\frac{1}{\sqrt{F+G}}\frac{\partial\zeta}{\partial u}du =
i\int^{\zeta(I,q_{2})}\frac{1}{\sqrt{F+G}}d\zeta.
$$
This can be written as a Jacobi elliptic integral; cf. \cite{Co, p.
400}.  Set  $z^2=(\zeta-\alpha)/(\zeta-\beta)$  where  $\alpha,\beta$
are the nonzero roots of  $F+G$;  then the last integral becomes
$$
\Theta_2 = \frac{i}{\sqrt{\alpha}}\int^{z(q_{2},I)}
\frac{dz}{\sqrt{(1-z^2)(1-k^2z^2)}},\qquad k^2=\beta/\alpha.\tag7.17
$$
Note that  $\alpha\beta = I_1I_2I_3$  is real and positive on  $SU(3)$.

The other angle variables are also elliptic integrals.
Straightforward calculation gives
$$\align
\Theta_1 &= \frac{1}{I_1}q_1 + i\int^{\zeta} \frac{1}{\sqrt{F+G}}
\{-1 + \frac{I_0+2I_3}{\zeta-I_1} + \frac{I_3(I_0+2)}{\zeta-I_1I_3}\}d\zeta;\\
\Theta_3 &= \frac{1}{I_3}q_3 + i\int^{\zeta}\frac{1}{\sqrt{F+G}} \{-1
+ \frac{I_0+2I_1}{\zeta-I_3} + \frac{I_1(I_0+2)}{\zeta-I_1I_3}\}
\endalign
$$
The first of the three integrals in each line are equal to
$-\Theta_2$  and the remaining two transform, under the same change
of variables as in (7.17), into sums of Jacobi elliptic integrals of
the first and third kinds.

To complete our discussion of  $SU(3)$  we consider the integration
of the flow with Hamiltonian  $I_2$,  in the original coordinate
system.  The Poisson bracket determined by the foliation of the form
$i\Omega$  differs from the Poisson bracket (5.11) by a factor  $-i$.
Thus the flow on  $SU(3)$  is given by
$$
\dot f = -i(\log I_2,f) = -i(\log a_{22}A_{22},f), \quad f\in
C^{\infty}(SL(n)).\tag7.18
$$
\enddemo

\proclaim{Theorem 7.3}  On  $SL(3,\bold C)$,  let
$$
\rho = (I_1I_2I_3)^{\frac12}, \quad \omega =
\frac{1}{2i}\log(a_{11}a_{22}a_{33}/A_{11}A_{22}A_{33}).
$$
These functions and the functions  $a_{jk}A_{jk}$  and
$\frac{1}{2i}\log(a_{jk}/A_{jk})$  are real on  $SU(3)$.  Under the
flow (7.17)  $\omega$  evolves according to the pendulum equation
$$
\ddot{\omega} = -2\rho\sin\omega.\tag7.19
$$
Moreover
\roster
\item"(7.20)"  each\quad  $a_{jk}A_{jk}$\quad  is an algebraic function of
$\cos\omega$  and the  $I_j$;
\item"(7.21)"  each time derivative of
$[\frac{1}{2i}\log(a_{jk}/A_{jk})]$  is an algebraic function of
$\cos\omega$  and the  $I_j$.
\endroster
\endproclaim

\demo{Proof}  A direct but somewhat tedious calculation using (7.18),
(5.11), and various identities for cofactors gives
$$\align
\ddot{\omega} &= \frac{1}{2i}\log(a_{22}/A_{22})\ddot{\hphantom{x}} =
i(a_{11}a_{22}a_{33} - A_{11}A_{22}A_{33})\\
&= -2\rho\sin\omega.
\endalign
$$
To obtain (7.20) we use the identities
$$
a_{j1}A_{k1} + a_{j2}A_{k2} + a_{j3}A_{k3} = \delta_{jk} =
a_{1j}A_{1k} + a_{2j}A_{2k} + a_{3j}A_{3k}
$$
which come from (7.05) to obtain
$$
a_{jk}A_{jk} + a_{kj}A_{kj} = 1 - I_j-I_k+I_{\ell},
$$
for distinct  $j,k,\ell$.  Also
$$\align
a_{jk}A_{jk}a_{kj}A_{kj} &= a_{jk}a_{kj}A_{jk}A_{kj} =
(a_{jj}a_{kk}-A_{\ell})(A_{jj}A_{kk}-a_{\ell\ell}) \\
&= I_jI_k+I_{\ell}-2\rho\cos\omega.
\endalign
$$
Therefore  $a_{jk}A_{jk}$  and  $a_{kj}A_{kj}$  are the roots of a
quadratic equation with coefficients which are polynomials in the
$I_j,\rho = (I_1I_2I_3)^{\frac12}$,  and  $\cos\omega$.

Finally, we consider  $\log(a_{12}/A_{12})$.  Another direct
calculation gives
$$\align
\frac{d}{dt}[\frac{1}{2i}\log(\frac{a_{12}}{A_{12}})] &= -\frac14
I_2\{\frac{a_{11}A_{21}}{a_{12}A_{22}} +
\frac{A_{11}a_{21}}{A_{12}a_{22}}\} \\
&= -\frac14\{a_{11}a_{22}(A_{11}A_{22}-a_{33}) +
A_{11}A_{22}(a_{11}a_{22}-A_{33})\} \\
&= \frac14\{\frac{A_{11}A_{22}A_{33} + a_{11}a_{22}a_{33} -
2I_1I_2}{a_{12}A_{12}}\} \\
&= \frac{1}{2a_{12}A_{21}}(\rho\cos\omega - I_1I_2).
\endalign
$$
The calculations for other such terms are similar, and  (7.21) follows,
using (7.20).

Note that on $SU(3)$,  the functions considered above are
$a_{jk}A_{jk} = \vert a_{jk}\vert^2$  and \linebreak
$\frac{1}{2i}\log(a_{jk}/A_{jk}) = \operatorname{arg}~ a_{jk}$. 
\enddemo

\head 8.  General nondegenerate $J$ \endhead

In this section we discuss the case of a spectral problem (3.2)
whose characteristic matrix  $J$  has  $n$  distinct eigenvalues but is
not otherwise constrained.  Thus we may assume
$$
J=\operatorname{diag}(i\lambda_1,i\lambda_2,\dots,i\lambda_n),\qquad
\lambda_j\text{'s \quad distinct}.\tag8.1
$$
The corresponding space of potentials  $P$, the $2$-form  $\Omega_p$,
and the associated Poisson bracket  $\{~,~\}_p$  are defined as in
section 2.  The (continuous) scattering data  $s$  or  $(v_+,v_-)$
which correspond to a potential  $q$ in  $P$  is a matrix-valued
function or pair of such functions, defined on the set
$$
\Sigma = \{\xi\in\bold C : Re(i\lambda_j\xi) =
Re(i\lambda_k\xi),\quad\text{some}\quad j\neq k\}.\tag8.2
$$
This set is a union of lines through the origin; we will consider it
as a union of rays from the origin and orient each ray from  $0$  to
$\infty$.  We will describe the spaces  $SD = \{s\}$  and  $SD' =
\{(v_+,v_-)\}$  in more detail below.  The analogue of Theorem 3.1
carries over to this more general setting; \cite{BC1}.  Therefore the
$2$-form and Poisson bracket can be carried over to a form
$\Omega_S$  and Poisson bracket  $\{~,~\}_S$  on scattering data.  In
this section we prove the existence of Darboux coordinates and the
complete integrability of the linear flows. 

\proclaim{Theorem 8.1}  There are functions  $p_{\nu},q_{\nu},1\leq
\nu\leq (n^2-n)/2$,  holomorphic on a dense open algebraic subset of
$SL(n)$, which have the following properties.  Let
$\Sigma_1,\dots,\Sigma_m$ be the rays of  $\Sigma$.  There is an
assignment of rays  $\nu\mapsto\Sigma_{k(\nu)}$  such that 
\roster
\item"(8.3)"  $s$  in  $SD$  is uniquely determined by the values
$$
\{p_{\nu}(s(\xi), q_{\nu}(s(\xi)),\quad  \xi\in\Sigma_{k(\nu)},\quad  1
\leq\nu\leq (n^2-n)/2\}.
$$
\item"(8.4)"  $\Omega_S = \sum\limits_{\nu}\int_{\Sigma_{k(\nu)}}\delta
p_{\nu}\wedge\delta q_{\nu}$.
\item"(8.5)"  $\{p_{\mu,\xi},p_{\nu,\eta}\}_S =
\{q_{\mu,\xi},q_{\nu,\eta}\}_S = 0$; \quad
$\{p_{\nu,\xi},q_{\nu,\eta}\}_S = \delta_{\mu\nu}\delta(\xi-\eta), \quad
\xi,\eta\in\Sigma_{k(\nu)}$. 
\item"(8.6)"  For any  $f$  in  $C^{\infty}(SL(n))$,  the
distributions   $\{f_{\xi},p_{\nu,\eta}\}$  have support at  $\xi =
\eta\in\Sigma_{k(\nu)}$,  all  $\nu$;  as before,  $f_{\xi}(s) =
f(s(\xi))$. 
\endroster
\endproclaim

\proclaim{Theorem 8.2}  The functions  $p_{\nu}$  of Theorem 8.1 may
be chosen so that for each traceless diagonal matrix  $\mu$, the
Hamiltonian for the flows (2.8), (3.20) is a linear combination of
the functionals 
$$
\int_{\Sigma_{k(\nu)}}\xi^kp_{\nu}(s(\xi))d\xi.
$$
\endproclaim

The machinery needed for the proofs of Theorems 8.1 and 8.2 has
already been developed in sections 3 and 4.  To show how it applies,
we must describe the space of scattering data.  Assume  $\int\Vert
q(x)\Vert dx<1$  and consider the spectral problem (3.1), (3.2).
Again there is a unique solution  $\psi(\cdot,z)$.  This solution is
holomorphic with respect to  $z, z\in\bold C\backslash\Sigma$, and
its boundary values satisfy 
$$
\psi(x,(1+i0)\xi) = \psi(x,(1-i0)\xi)v(\xi),\quad  \xi\in\Sigma\backslash
0.\tag8.7
$$
Given  $\xi$  in  $\Sigma\backslash 0$, let  $\Pi_{\xi}:M_n\to M_n$
be the projection defined by
$$
(\Pi_{\xi}a)_{jk} = \cases a_{jk} &\quad\text{if}\quad
Re(i\lambda_j\xi) = Re(i\lambda_k\xi),\\ 0, &\quad\text{otherwise}.\endcases
$$
Then the limit
$$
\lim\limits_{x\to +\infty} \Pi_{\xi}\psi(x,\xi)e^{-x\xi J} = s(\xi).\tag8.8
$$
exists on  $\Sigma\backslash 0$.  There are factorizations
$$
v(\xi) = v_-(\xi)^{-1}v_+(\xi) = s_-(\xi)^{-1}s_+(\xi);\quad s_{\pm}(\xi)
= s(\xi)v_{\pm}(\xi).\tag8.9
$$
The factors in (8.9) are characterized by the following conditions:
$$\align
\Pi_{\xi}v_{\pm}(\xi) &= v_{\pm}(\xi); \quad \Pi_{\xi}s_{\pm}(\xi) =~~
s_{\pm}(\xi); \quad \operatorname{diag}~v_{\pm} = 1;\tag8.10\\
(v_{\pm}(\xi))_{jk} &= (s_{\mp}(\xi))_{jk} = 0 \\
\text{if} &\quad Re~i(\lambda_j-\lambda_k)w > 0\quad\text{for}
\quad w=(1+i\varepsilon)\xi,\quad\text{small}\quad \varepsilon
>0.
\endalign
$$
There are further conditions on  $s$  and  $v$  as functions of
$\xi\in\Sigma$,  which we do not need to cite here; see [BC1] for a
complete discussion of these conditions and of the algebraic facts we
are using in this section.

Proposition 3.2 carries over, in the following form \cite{BC3}:
$$
\Omega_S = \frac{1}{2\pi i}\int_{\Sigma}tr(v_-(\delta
v)v^{-1}_+\wedge s^{-1}\delta s)dz.\tag8.11
$$

The strategy for proving Theorems 8.1 and 8.2 is the same as for
proving Theorems 3.7 and 3.8: pointwise analysis of the form
$\Omega_S$.

Suppose  $\xi$  is in  $\Sigma\backslash 0$.  After conjugation by a
permutation matrix (which depends on the ray of  $\Sigma$  containing
$\xi$), we may assume that
$$\align
Re~i\lambda_1\xi &= \dots = Re~i\lambda_{d_{1}}\xi >
Re~i\lambda_{d_{1}+1}\xi = \dots\tag8.12\\
 &= Re~i\lambda_{d_{2}}\xi > Re~i\lambda_{d_{2}+1}\xi \dots~~.
\endalign
$$
Then  $v(\xi)$  and  $s(\xi)$  have the block diagonal form
$$\bmatrix
A_1 & 0   &\dotso & 0\\
0   & A_2 &\dotso & 0\\
&\dotso\\
0 & 0 &\dotso & A_s\endbmatrix ,\qquad A_j\in SL(d_j).\tag8.13
$$
Moreover   $v_-(\xi)$  and  $s_+(\xi)$  have this form and are upper
triangular, while  $v_+(\xi)$  and  $s_-(\xi)$  have this form and
are lower triangular.  With this normalization, the pointwise
$2$-form being integrated over the ray containing  $\xi$  is (a
constant multiple of) the sum of the forms as in Theorem 4.1 for the
matrix groups  $GL(d_1), GL(d_2),\dots$.  There is an analogous
decomposition of the Poisson bracket, which is computed as in section
3.

It follows from these considerations that the (trivial) extension
of the results in sections 4 and 5 to block diagonal matrix groups
yields the desired functions  $p_{\nu}, q_{\nu}$  of Theorems 8.1 and
8.2.

\vskip .3in

\centerline{\bf References}
\vskip .25in

\roster

\item"\bf [AKNS]\rm"  M. J. Ablowitz, D. J. Kaup, A. C. Newell, and H.
Segur, The inverse scattering transform--Fourier analysis for nonlinear
problems, Stud.  Appl. Math. 53 (l974), 249-315.
\vskip .10in

\item"\bf [A]\rm"  V. I. Arnold, Mathematical Methods of Classical
Mechanics, Springer, New York, 1989.

\item"\bf [BY]\rm"    D. Bar Yaacov, Analytic properties of scattering and
inverse scattering for first order systems, Dissertation, Yale l985.
\vskip .10in

\item"\bf [BC1]\rm"   R. Beals and R. R. Coifman, Scattering and inverse
scattering for first order systems, Comm. Pure Appl. Math. 37 (1984),
39-90. 
\vskip .10in

\item"\bf [BC2]\rm"   R. Beals and R. R. Coifman, Inverse scattering and
evolution equations, Comm. Pure Appl. Math. 38 (l985), 29-42.
\vskip .10in

\item"\bf [BC3]\rm"   R. Beals and R. R. Coifman, Linear spectral problems,
nonlinear equations, and the  $\bar{\partial}$-method, Inverse Problems 5
(l989), 87-130.
\vskip .10in

\item"\bf [Ca]\rm"    P. J. Caudrey, The inverse problem for a general
$n\times n$  spectral equation, Physica D6 (l982), 51-66.
\vskip .10in

\item"\bf [Co]\rm"    E. T. Copson, Theory of Functions of a Complex Variable,
Oxford, 1935.
\vskip .10in

\item"\bf [Dr]\rm"   V. G. Drinfeld, "Quantum groups", in Proc. International
Congress of Mathematicians, Berkeley, l986.
\vskip .10in

\item"\bf [Ga]\rm"    C. S. Gardner, Korteweg-de Vries equation and
generalizations, IV.  The Korteweg-de Vries equation as a Hamiltonian
system, J. Math.  Physics 12 (1971), 1548-1551.
\vskip .10in

\item"\bf [Ge]\rm"    V. S. Gerdzhikov, On the spectral theory of the
integro-differential operator  $\Lambda$  generating nonlinear evolution 
equations, Lett. Math.  Physics 6 (1982), 315-323. 
\vskip .10in

\item"\bf [Ka]\rm"    D. J. Kaup, The three-wave interaction -- a 
non-dispersive phenomenon, Studies Appl. Math. 55 (l976), 9-44.
\vskip .10in

\item"\bf [KD]\rm"  B. G. Konopelchenko and V. G. Dubrovsky, Hierarchy of
Poisson brackets for elements of a scattering matrix, Lett. Math.
Phys 8 (1984), 273.
\vskip .10in

\item"\bf [La]\rm"  P. D. Lax, Almost periodic solutions of the  $KdV$
equation, SIAM Review 18 (1976), 351-375.
\vskip .10in

\item"\bf [Lu]\rm"    J-h. Lu, personal communication.
\vskip .10in

\item"\bf [Ma]\rm"   S. V. Manakov, An example of a completely integrable
nonlinear wave field with nontrivial dynamics (Lee Model), Teor. Mat.
Phys. 28 (1976), 172-179.
\vskip .10in

\item"\bf [Ne]\rm"    A. C. Newell, The general structure of
integrable evolution equations, Proc. R. Soc. Lond. A 365 (1979),
283-311.
\vskip .10in

\item"\bf [Sa]\rm"    D. H. Sattinger, Hamiltonian hierarchies on
semisimple Lie algebras, Stud. Appl. Math. 72 (l985), 65-86.
\vskip .10in

\item"\bf [Sh]\rm"    A. B. Shabat, An inverse scattering problem, Diff.
Uravn. 15 (1978), l824-1834; Diff. Equ. 15 (1980), 1299-1307.
\vskip .10in

\item"\bf [Sk]\rm"    E. K. Sklyanin, Quantum variant of the method
of the inverse scattering problem, Zap. Nauchn. Sem. Leningrad.
Otdel. Mat. Inst. Steklov 95 (1980), 55-128, 161.
\vskip .10in

\item"\bf [Wh]\rm"  E. T. Whittaker, Analytical Mechanics, Dover, New York,
l944.  
\vskip .10in

\item"\bf [ZF]\rm"  V. E. Zakharov and L. D. Faddeev, Korteweg-de Vries
equation: a completely integrable Hamiltonian system, Funct. Anal.
Appl. 5 (1971), 280-287.
\vskip .10in

\item"\bf [ZM1]\rm"  V. E. Zakharov and S. V. Manakov, On resonant
interaction of wave packets in nonlinear media, Zh\'ETF Pis.Red. 18
(1973), 413.
\vskip .10in

\item"\bf [ZM2]\rm"  V. E. Zakharov and S. V. Manakov, On the
complete integrability of the nonlinear Schr\"odinger equation, Teor.
Mat. Fyz. 19 (1974), 332-343.
\vskip .10in

\item"\bf [ZS1]\rm"  V. E. Zakharov and A. B. Shabat, Exact theory of
two-dimensional self-focussing and one-dimensional self-modulation of
waves in nonlinear media, Soviet Physics JETP 34 (l972), 62-69.
\vskip .10in

\item"\bf [ZS2]\rm"  V. E. Zakharov and A. B. Shabat, A scheme for integrating
nonlinear equations of mathematical physics by the method of the
inverse scattering transform, I, Funct. Anal. Appl. 8 (l974),
226-235.

\endroster

\end